\documentclass[12pt,a4paper]{amsart}
\usepackage{amscd,amsmath,amssymb,amsfonts,times}
\usepackage{mathrsfs}
\usepackage[all]{xy}
\numberwithin{equation}{subsection}
\theoremstyle{plain}
    \newtheorem{thm}{Theorem}[subsection]
    \newtheorem{lem}[thm]{Lemma}
    
    \newtheorem{sublem}[thm]{Sublemma}
    \newtheorem{prop}[thm]{Proposition}
    \newtheorem{cor}[thm]{Corollary}

\theoremstyle{definition}
    \newtheorem{defn}[thm]{Definition}

\theoremstyle{remark}
    \newtheorem{rem}[thm]{Remark}
\newenvironment{pf}{\par\smallskip\noindent\emph{Proof.}}{\qed\par\medskip}
\newenvironment{pf*}[1]{\par\smallskip\noindent\emph{#1.}}{\qed\par\medskip}

\begin{document}
\title[zero-cycles and Brauer groups]{zero-cycles on varieties over $\boldsymbol{p}$-adic fields \endgraf and Brauer groups}
\author[S. Saito and K. Sato]{Shuji Saito and Kanetomo Sato}
\date{December 19, 2013}
\address[Shuji Saito]{\;Department of Mathematics, Tokyo Institute of Technology \endgraf
 2-12-1 Oh-okayama, Meguro-ku, Tokyo 152-8551, JAPAN}
\email{sshuji@msb.biglobe.ne.jp}
\address[Kanetomo Sato]{\;Department of Mathematics, Chuo University \endgraf
 1-13-27 Kasuga, Bunkyo-ku, Tokyo 112-8551, JAPAN}
\email{kanetomo@math.nagoya-u.ac.jp}
\begin{abstract} In this paper, we study the Brauer-Manin pairing of smooth proper varieties over a $p$-adic field, and determine the $p$-adic part of the image of the induced cycle map. We also compute {\it A}$_0$ of a potentially rational surface which splits over a wildly ramified extension.
\par \bigskip  \bigskip \noindent
{\sc R\'esum\'e.}  Dans cet article, nous \'etudions l'accouplement de Brauer-Manin des vari\'et\'es propres et lisses sur un corps $p$-adique, et d\'eterminons la partie $p$-adique de l'image de l'application cycle induite. Nous calculons aussi {\it A}$_0$ d'une surface potentiellement rationelle qui se scinde sur une extension sauvagement ramifie\'e.
\end{abstract}
\maketitle
\setlength{\baselineskip}{15pt}

%
%

%
%
%
\def\A{\text{\it A}}        
\def\ab{\text{ab}}      
\def\b{\text b}          
\def\Br{\text{\rm Br}}      
\def\br{\text{\rm br}}      
\def\cd{\text{\rm cd}}      
\def\codim{\text{\rm codim}}
\def\ch{\text{\rm ch}}      
\def\CH{\text{\rm CH}}      
\def\cl{\text{\rm cl}}      
\def\Coker{\text{\rm Coker}}
\def\Ker{\text{\rm Ker}}    
\def\cont{\text{\rm cont}}  
\def\Cor{\text{\rm Cores}}    
\def\div{\text{\rm div}}    
\def\dlog{d{\text{\rm log}}}  
\def\et{\text{\rm \'et}}    
\def\F{\text{\rm F}}        
\def\Fil{\text{\rm Fil}}    
\def\fin{\text{\rm fin}}    
\def\Frac{\text{Frac}}  
\def\G{\text{G}}        
\def\Gal{\text{Gal}}    
\def\gcd{\text{G.C.D.}} 
\def\gr{\text{gr}}      
\def\gys{\text{Gys}}    
\def\h{{\hspace{.4pt}{\sf h}}}        
\def\H{H}        
\def\Hom{\text{\rm Hom}}    
\def\id{\text{\rm id}}     
\def\Image{\text{\rm Im}}   
\def\inv{\text{\rm inv}}    
\def\ip{[p^{-1}]}
\def\ker{\text{\rm Ker}}    
\def\log{\text{\rm log}}    
\def\ltor{{\ell\text{-}\text{\rm tors}}}
\def\mod{\;{\text{mod}}\;}    
\def\NF{N}       
\def\NS{{\text{\rm NS}}}      
\def\ord{{\text{\rm ord}}}    
\def\pDiv{p\text{-}\Div}    
\def\Pic{{\text{\rm Pic}}}    
\def\pro{{\text{\rm pro}\text{-}}} 
\def\Proj{{\text{\rm Proj}}}  
\def\ptor{{p\text{-}\text{\rm tors}}}
\def\red{{\text{\rm red}}}    
\def\res{{\text{\rm res}}}    
\def\Res{{\text{\rm Res}}}    
\def\sh{{\hspace{.8pt}\text{\rm sh}}}      
\def\sing{{\text{\rm sing}}}  
\def\Spec{{\text{\rm Spec}}}  
\def\st{{\text{\rm st}}}      
\def\sym{{\text{\rm sym}}}    
\def\tor{{\text{\rm tors}}}   
\def\Cotor{{\text{\rm Cotor}}}   
\def\tr{{\text{\rm tr}}}      
\def\tsym{{\text{tsym}}}  
\def\ur{{\text{\rm ur}}}      
\def\sp{{\text{\rm sp}}}      

\def\BB{{\mathscr B}}
\def\cC{{\mathscr C}}
\def\FF{{\mathscr F}}           
\def\cI{{\mathscr I}}
\def\cJ{{\mathscr J}}
\def\K{{\mathscr K}}            
\def\cK{{\mathscr K}}
\def\cL{{\mathscr L}}
\def\O{{\mathscr O}}            
\def\cO{{\mathscr O}}
\def\cS{{\mathscr S}}
\def\cU{{\mathscr U}}
\def\V{{\mathscr V}}
\def\WW{{\mathscr W}}
\def\ZZ{{\mathscr Z}}

\def\bA{{\mathbb A}}
\def\bF{{\mathbb F}}
\def\bP{{\mathbb P}}
\def\hH{{\mathbb{H}}}
\def\L{{\mathbb L}}
\def\N{{\mathbb N}}
\def\Q{{\mathbb Q}}
\def\R{{\mathbb R}}
\def\Z{{\mathbb Z}}
\def\rf{{\mathbb F}}

\def\fI{{\mathfrak I}}
\def\fJ{{\mathfrak J}}
\def\T{{\mathfrak T}}
\def\X{{\mathfrak X}}
\def\fm{{\mathfrak m}}
\def\fo{{\mathfrak o}}
\def\fp{{\mathfrak p}}
\def\OO{{\mathfrak o}}
\def\QQ{{\mathfrak C}}
\def\p{{\mathfrak p}}

\def\btil{\wt{\phantom{a}}}
\def\ctil{\hspace{-4pt}\wt{\phantom{a}}}
\def\dtil{\hspace{-10pt}\wt{\phantom{a}}\hspace{3pt}}

%
\def\Lam{\varLambda}
\def\lam{\lambda}
\def\Te{\Theta}
\def\te{\theta}
\def\vare{\varepsilon}
\def\vG{\Gamma}
\def\ze{\zeta}
\def\k{\kappa}
%
%
%
%
\def\ra{\rightarrow}
\def\lra{\longrightarrow}
\def\Lra{\Longrightarrow}
\def\la{\leftarrow}
\def\lla{\longleftarrow}
\def\Lla{\Longleftarrow}
\def\da{\downarrow}
\def\hra{\hookrightarrow}
\def\lmt{\longmapsto}
\def\sm{\setminus}
\def\wt#1{\widetilde{#1}}
\def\wh#1{\widehat{#1}}
\def\spt{\sptilde}
\def\ol#1{\overline{#1}}
\def\ul#1{\underline{#1}}
\def\us#1#2{\underset{#1}{#2}}
\def\os#1#2{\overset{#1}{#2}}
\def\lim#1{\us{#1}{\varinjlim}}
\def\bs#1{\boldsymbol{#1}}
\def\rmapo#1{\overset{#1}{\longrightarrow}}
%
%
%
%
\def\a{\alpha}
\def\b{\beta}
\def\c{\psi}
\def\d{\delta}
\def\f{\varphi}
\def\n{\nu}
\def\o{\rho}
\def\q{\gamma}
\def\r{\theta}
\def\t{\tau}
\def\w{\omega}
\def\W{\Omega}
\def\aa{{\text{a}}}
\def\ab{\overline{a}}
\def\an{\alpha_\n}
\def\Ah{A_{\eta}}
\def\Ahx{A_{\eta_x}}
\def\Ahxb{A_{\eta_{\ol x}}}
\def\bb{\overline{b}}
\def\bn{\beta_\n}
\def\cb{\overline{c}}
\def\hx{\eta_x}
\def\hxb{{\eta_{\ol x}}}
\def\ub{\overline{u}}
\def\Kh{K_{\eta}}
\def\Khx{K_{\eta_x}}
\def\Wo#1{\Omega^1_{#1}}
\def\Wt#1{\Omega^2_{#1}}
\def\Wtlog#1{\Omega^2_{#1,\log}}
\def\Wwolog#1#2{W_{#1}\omega^1_{#2,\log}}
\def\WWolog#1#2{W_{#1}\Omega^1_{#2,\log}}

\def\Zo#1{{\mathscr Z}^1_{#1}}
\def\Bo#1{B^1_{#1}}
\def\WZo#1{\Omega^1_{#1}/{\mathscr Z}^1_{#1}}
\def\Gm{{\mathbb G}_{\hspace{-1pt}\text{m}}}
\def\Ga{{\mathbb G}_{\hspace{-1pt}\text{a}}}
%
%
%
\def\mwitt#1#2#3{W_{\hspace{-2pt}#2}{\hspace{1pt}}\omega_{#1}^{#3}}
\def\witt#1#2#3{W_{\hspace{-2pt}#2}{\hspace{1pt}}\Omega_{#1}^{#3}}
\def\mlogwitt#1#2#3{W_{\hspace{-2pt}#2}{\hspace{1pt}}\omega_{{#1},{\log}}^{#3}}
\def\logwitt#1#2#3{W_{\hspace{-2pt}#2}{\hspace{1pt}}\Omega_{{#1},{\log}}^{#3}}
%
%
\def\bR{{\mathbb R}}
\def\bZ{{\mathbb Z}}
\def\bQ{{\mathbb Q}}
\def\qz{{\bQ}/{\bZ}}
\def\qzb{\Q/\Z'}
\def\qzl{{\bQ_\ell}/{\bZ_\ell}}
\def\qzp{{\bQ_p}/{\bZ_p}}
\def\zl{{\bZ}_\ell}
\def\zp{{\bZ}_p}
\def\pnuz{\bZ/p^r\bZ}

\def\knr{k^{\ur}}
\def\Xnr{X^{\ur}}
\def\Ynr{\ol Y}
\def\gfiber{X}
\def\cX{{\mathscr X}}
\def\scX{\hspace{-1.3pt}\cX\hspace{-1.3pt}}
\def\scXp{\hspace{-1.3pt}\cX'\hspace{-1.3pt}}
\def\sscX{\hspace{-1.7pt}\cX\hspace{-1.7pt}}
\def\isom{\hspace{9pt}{}^\sim\hspace{-16.5pt}\lra}
\def\lisom{\hspace{9pt}{}^\sim\hspace{-17.5pt}\lla}
\def\psy{\hspace{.5pt}[\hspace{-1.57pt}[y]\hspace{-1.58pt}]}
\def\psxy{\hspace{.5pt}[\hspace{-1.57pt}[x,y]\hspace{-1.58pt}]}

\def\cP{{\mathscr P}}
\def\cQ{{\mathscr Q}}
\def\cR{{\mathscr R}}
\def\indlim#1{\underset{#1}{\varinjlim}\ }
\def\projlim#1{\underset{#1}{\varprojlim}\ }

\section{Introduction}\label{sect1}
\medskip
Let $k$ be a $p$-adic local field, and let $\gfiber$ be a proper smooth geometrically integral variety over $k$.
Let $\CH_0(\gfiber)$ be the Chow group of $0$-cycles on $\gfiber$ modulo rational equivalence.
An important tool to study $\CH_0(\gfiber)$ is the natural pairing due to Manin \cite{M}
\begin{align}
\tag{M}\label{eq01} \CH_0(\gfiber) \times \Br(\gfiber) \lra  \Q/\Z,
\end{align}
where $\Br(\gfiber)$ denotes the Grothendieck-Brauer group $\H^2_{\et}(\gfiber,\Gm)$. When $\dim (\gfiber) = 1$, using the Tate duality theorem for abelian varieties over $p$-adic local fields, Lichtenbaum \cite{Li2} proved that \eqref{eq01} is non-degenerate and induces an isomorphism
\begin{equation}\tag{L}\label{eq011}
 \A_0(\gfiber) \isom \Hom(\Br(\gfiber)/\Br(k), \Q/\Z).
\end{equation}
Here $\Br(\gfiber)/\Br(k)$ denotes the cokernel of the natural map $\Br(k) \ra \Br(\gfiber)$, and $\A_0(\gfiber)$ denotes the subgroup of $\CH_0(\gfiber)$ generated by $0$-cycles of degree $0$.
An interesting question is as to whether the pairing \eqref{eq01} is non-degenerate when $\dim(\gfiber) \ge 2$. See \cite{PS} for surfaces with non-zero left kernel. See \cite{Y2} for varieties with trivial left kernel. In this paper, we are concerned with the right kernel of \eqref{eq01} in the higher-dimensional case.
\subsection{}\label{sect1-1}
We assume that $\gfiber$ has a regular model $\cX$ which is proper flat of finite type over the integer ring $\OO_k$ of $k$. It is easy to see that the pairing \eqref{eq01} induces homomorphisms
\begin{align}
\label{eq02} \CH_0(\gfiber) & \lra \Hom(\Br(\gfiber)/\Br(\cX), \Q/\Z), \\
\label{eq03} \A_0(\gfiber) & \lra \Hom(\Br(\gfiber)/\Br(k)+\Br(\cX), \Q/\Z),
\end{align}
where $\Br(\gfiber)/\Br(k)+\Br(\cX)$ denotes the quotient of $\Br(\gfiber)$ by the image of $\Br(k)\oplus \Br(\cX)$.
If $\dim(\gfiber)=1$, then $\Br(\cX)$ is zero, and the map \eqref{eq03} is the same as \eqref{eq011} (cf.\ \cite{CTOP} 1.7\,(c)). Our main result is the following:
\addtocounter{thm}{2}
\begin{thm}\label{cor0-3}
Assume that the purity of Brauer groups holds for $\cX$ {\rm(}see Definition {\rm\ref{def:purity}} below{\rm)}.
Then{\rm:}
\begin{enumerate}
\item[(1)]
The right kernel of the pairing \eqref{eq01} is exactly $\Br(\cX)$,
  that is, the map \eqref{eq02} has dense image
   with respect to the natural pro-finite topology on
     the right hand side.
\item[(2)] The map \eqref{eq03} is surjective.
\end{enumerate}
\end{thm}
\noindent
Restricted to the prime-to-$p$ part, the assertion (1) is due to Colliot-Th\'el\`ene and Saito \cite{cts}. The assertion (2) gives an affirmative answer to \cite{CT2} Conjecture 1.4\,(c), assuming the purity of Brauer groups, which holds if $\dim(\cX) \le 3$ or if $\cX$ has good or semistable reduction (cf.\ Remark \ref{rem:purity} below).
Roughly speaking, Theorem \ref{cor0-3}\,(1) asserts that if an element $\omega \in \Br(\gfiber)$ ramifies along the closed fiber of $\cX$, then there exists a closed point $x \in \gfiber$ for which the specialization of $\omega$ is non-zero in $\Br(x)$. We will in fact prove the following stronger result on the ramification of Brauer groups:
\begin{thm}[{{\bf Corollary \ref{cor:key}}}]\label{thm0-2}
Let $\cU$ be either $\cX$ itself or its henselization at a closed point.
Put $U:=\cU[p^{-1}]$ and assume that the purity of Brauer groups holds for $\cU$. If $\cX$ is henselian local, then assume further that all irreducible components of the divisor on $\cU$ defined by the radical of $(p)$ are regular. Then the kernel of the map
\[ \psi_x : \Br(U) \lra \prod_{v\in U_0} \ \qz\; , \quad \omega \mapsto (\inv_v(\omega\vert_v))_{v\in U_0} \]
agrees with $\Br(\cU)$, where $U_0$ denotes the set of closed points on $U$.
\end{thm}
\noindent
The prime-to-$p$ part of Theorem \ref{thm0-2} has been proved in \cite{cts}. We will prove the $p$-primary part of this result using Kerz's id\`ele class group \cite{Ke}. Our method of the proof gives also an alternative proof of the prime-to-part in \cite{cts}.

\subsection{}\label{sect1-1'}
As an application of Theorem \ref{cor0-3}\,(2), we give an explicit calculation of $\A_0(\gfiber)$ for a potentially rational suface $\gfiber/k$, a proper smooth geometrically connected surface $\gfiber$ over $k$ such that $\gfiber \otimes_k k'$ is rational for some finite extension $k'/k$. For such a surface $X$, the map \eqref{eq03} has been known to be injective (see Proposition \ref{prop.ex1} below), and hence bijective by Theorem \ref{cor0-3}\,(2). On the other hand, for such a surface $X$, we have
\[ \Br(\gfiber)/\Br(k) \simeq H^1_{\Gal}(G_k,\NS(\ol{\gfiber})), \]
where $\NS(\ol{\gfiber})$ denotes the N\'eron-Severi group of $\ol{\gfiber}:=\gfiber \otimes_k \ol k$, and $G_k$ denotes the absolute Galois group of $k$. Thus knowing the $G_k$-module structure of $\NS(\ol{\gfiber})$, we can compute $\A_0(\gfiber)$ by determining which element of $\Br(\gfiber)$ are unramified along the closed fiber of $\cX$.
For example, consider a cubic surface for $a\in k^\times$ \[ \gfiber\, : \, T_0^3 + T_1^3 + T_2^3 + a T_3^3 =0 \quad
\hbox{ in \;\; $\bP^3_k=\Proj (k[T_0,T_1,T_2,T_3])$}. \]
If $a$ is a cube in $k$, then $X$ is isomorphic to the blow-up of $\bP^2_k$ at six $k$-valued points in the general position (Shafarevich) and we have $\A_0(\gfiber)=0$. We will prove the following result, which is an extention of results in \cite{cts} Example 2.8.
\begin{thm}[{{\bf Theorem \ref{thm.cubic}}}]\label{thm0-3}
Assume that $\ord_k(a) \equiv 1 \mod (3)$ and that $k$ contains a primitive cubic root of unity. Then we have \[ \A_0(\gfiber) \simeq  (\bZ/3)^2. \]
\end{thm}
\noindent
In his paper \cite{Da}, Dalawat provided a method to compute $\A_0(\gfiber)$ for a potentially rational surface $\gfiber$, which works under the assumption that the action of $G_k$ on $\NS(\ol{\gfiber})$ is unramified. Theorem \ref{cor0-3} provides a new method to compute $\A_0(\gfiber)$, which does not requires Dalawat's assumption. Note that $p$ may be $3$ in Theorem \ref{thm0-3}, so that the action of $G_k$ on $\NS(\ol{\gfiber})$ may ramify even wildly.
\subsection{}\label{sect1-2}
Let $\OO_k$ be as before, and let $\cX$ be a regular scheme which is proper flat of finite type over $\OO_k$.
Assume that $\cX$ has {\it good or semistable reduction} over $\OO_k$.
Let $d$ be the absolute dimension of $\cX$, and let $r$ be a positive integer.
In \cite{SS}, we proved that the cycle class map
\[ \varrho_m^{d-1} : \CH^{d-1}(\cX)/m  \lra \H^{2d-2}_{\et}(\cX,\mu_m^{\otimes d-1}) \]
 is bijective for any positive integer $m$ prime to $p$. Here $\mu_m$ denotes the \'etale sheaf of $m$-th roots of unity.
As a new tool to study $\CH^{d-1}(\cX)$, we introduce the $p$-adic cycle class map defined in \cite{Sat2} Corollary 6.1.4:
\[ \varrho_{p^r}^{d-1} : \CH^{d-1}(\cX)/p^r  \lra \H^{2d-2}_{\et}(\cX,\T_r(d-1)). \]
Here $\T_r(n)=\T_r(n)_{\scX}$ denotes the \'etale Tate twist with $\Z/p^r\Z$-coefficients \cite{Sat2} (see also \cite{Sch} \S7),
  which is an object of $D^b(\cX,\pnuz)$, the derived category of bounded complexes of \'etale $\Z/p^r\Z$-sheaves on $\cX$.
This object $\T_r(n)$ plays the role of $\mu_m^{\otimes n}$, and we expect that $\T_r(n)$ agrees with $\Z(n)^{\et} \otimes^{\mathbb L} \Z/p^r\Z$, where $\Z(n)^{\et}$ denotes the conjectural \'etale motivic complex of Beilinson-Lichtenbaum
  (\cite{Be}, \cite{Li}, \cite{Sat2} Conjecture 1.4.1\,(1)).
Concerning the map $\varrho_{p^r}^{d-1}$, we will prove the following result:
\begin{thm}\label{thm:surj}
The cycle class map $\varrho_{p^r}^{d-1}$ is surjective.
\end{thm}
\noindent
We have nothing to say about the injectivity of $\varrho_{p^r}^{d-1}$ in this paper (compare with \cite{Y1}).
A key to the proof of Theorem \ref{thm:surj} is the non-degeneracy of a canonical pairing of finite $\Z/p^r\Z$-modules
\[ \H^{2d-2}_{\et}(\cX,\T_r(d-1)) \times \H^3_{Y,\et}(\cX,\T_r(1)) \lra \Z/p^r\Z \]
proved in \cite{Sat2} Theorem 10.1.1. We explain an outline of the proof of Theorem \ref{thm:surj}.
Let $Y$, $U$, $A_x$ be as in Theorem \ref{thm0-2}. Let $\gfiber_0$ and $Y_0$ be the sets of all closed points on $\gfiber$ and $Y$, respectively, and let $\sp:\gfiber_0 \to Y_0$ be the specialization map of points. By the duality mentioned above, there is an isomorphism of finite groups
\[ \H^{2d-2}_{\et}(\cX,\T_r(d-1)) \isom \H^3_{Y,\et}(\cX,\T_r(1))^*, \]
where we put $M^*:=\Hom(M,\qz)$ for abelian group $M$. We will construct an injective map
\[\xymatrix{ \us{\phantom{1}}{\theta_{p^r} : \H^3_{Y,\et}(\cX,\T_r(1))} \; \ar@<-1pt>@{^{(}->}[r] & \displaystyle \prod_{x\in U_0} \ {}_{p^r}\Br(A_x[p^{-1}]) }\]
whose dual fits into a commutative diagram
\[\xymatrix{ \CH^{d-1}(\cX)/p^r \ar[r]^-{\varrho_{\ell^r}^{d-1}} &  \H^{2d-2}_{\et}(\cX,\T_r(d-1)) \ar[r] \ar@{}@<-.5mm>[r]^-{\sim} & \H^3_{Y,\et}(\cX,\T_r(1))^* \\
 \displaystyle \bigoplus_{x \in U_0} \ \bigoplus_{v\in \Spec(A_x[p^{-1}])_0} \ \Z/p^r\Z \ar@{->>}[rr]^-{(\psi_{p^r})^*} \ar[u] & & \displaystyle \bigoplus_{x\in U_0} \ \big({}_{p^r}\Br(A_x[p^{-1}])\big)^* \ar@{->>}[u]_{(\theta_{p^r})^*}. }\]
Here $\psi_{p^r}$ denotes the direct product of the $p^r$-torsion part of the map $\psi_x$ in Theorem \ref{thm0-2} for all $x \in U_0$, which is injective by Theorem \ref{thm0-2} and its dual $(\psi_{p^r})^*$ is surjective. Therefore Theorem \ref{thm:surj} will follow from this commutative diagram and the surjectivity of $(\theta_{p^r})^*$ and $(\psi_{p^r})^*$ (see \S\ref{sect6} for details).
\subsection{}
This paper is organized as follows. In \S\ref{sect4}, we will prove Theorem \ref{thm0-2} in a stronger form. In \S\ref{sect.cubic}, we compute $\A_0$ of cubic surfaces to prove Theorem \ref{thm0-3}. In \S\ref{sect5} and \S\ref{sect6}, we will prove Theorem \ref{cor0-3} and Theorem \ref{thm:surj}, respectively.
\par
\vspace{8pt}
\noindent
{\it Acknowledgements.}
\hspace{1pt}
The research for this article was partially supported by EPSRC grant and JSPS Core-to-Core Program.
The authors express their gratitude to Professor Ivan Fesenko for valuable comments and discussions 
and to The University of Nottingham for their great hospitality. Thanks are also due to Professor 
Jean-Louis Colliot-Th\'el\`ene for his lucid expository paper \cite{CTb} on the subject of this paper 
based on his Bourbaki talk, and for allowing them to include his note \cite{ct:letter} in this paper. 
They also thank Tetsuya Uematsu for several valuable comments on the arguments in \S\ref{sect.cubic}. Finally they 
express sincere gratitude to Professor Moritz Kerz for his valuable suggestion to simplify the proof of Theorem \ref{thm0-2}.
\section*{Notation}
\subsection{}
For an abelian group $M$ and a positive integer $n$,
     ${}_nM$ and $M/n$ denote the kernel and the cokernel of the map
       $M \os{\times n}{\lra} M$, respectively.
For a field $k$, $\ol k$ denotes a fixed separable closure, and
   $G_k$ denotes the absolute Galois group $\Gal(\ol k/k)$.
For a discrete $G_k$-module $M$, $\H^*(k,M)$ denote
   the Galois cohomology groups $\H^*_{\Gal}(G_k,M)$,
     which are the same as the \'etale cohomology groups of $\Spec(k)$
       with coefficients in the \'etale sheaf associated with $M$.
\subsection{}
Unless indicated otherwise, all cohomology groups of schemes are taken over the \'etale topology.
For a commutatitive ring $R$ with unity and a sheaf $\FF$ on $\Spec(R)_{\et}$,
   we often write $\H^*(R,\FF)$ for $\H^*(\Spec(R),\FF)$.
\subsection{}
For a scheme $X$, a sheaf $\FF$ on $X_{\et}$ and a point $x \in X$,
   we often write $\H^*_x(X,\FF)$ for $\H^*_x(\Spec(\O_{X,x}),\FF)$.
For a point $x \in X$, $\kappa(x)$ denotes its residue field.
We often write $X_0$ for the set of all closed points on $X$.
For a pure-dimensional scheme $X$ and a non-negative integer $q$,
    $X^q$ denotes the set of all points on $X$ of codimension $q$.
For an integer $n \geq 0$ and a noetherian excellent scheme $X$,
 $\CH_n(X)$ denotes the Chow group of algebraic cycles on $X$ of {\it dimension} $n$ modulo rational equivalence;
 if $X$ is regular of pure dimension $d$, we often write $\CH^n(X)$ for $\CH_{d-n}(X)$.
 
\section{Preliminaries}
In this section, we introduce some terminology and notions which will be useful throughout this paper.

\subsection{Purity of Brauer groups}

\begin{defn}\label{def:purity}
Let $X$ be a noetherian scheme.
\begin{enumerate}
\item[(1)]
For a closed immersion $\iota_Z:Z \ra X$ with $\codim_X(Z) \geq 2$,
   we say that {\it the purity of Brauer groups holds for the pair $(X,Z)$}, if $R^3\iota_Z^!\Gm,_X=0$.
\item[(2)]
We say that {\it the purity of Brauer groups holds for} $X$, if the purity of Brauer groups holds for
     any pair $(X,Z)$ with $\codim_X(Z) \geq 2$.
\end{enumerate}
\end{defn}
\begin{rem}\label{rem:purity}
There are some known cases on this purity problem:
\begin{enumerate}
\item[(1)]
For a noetherian regular scheme $X$ with $\dim(X) \leq 3$, the purity of Brauer groups holds for $X$ (\cite{gabber}).
\item[(2)]
For a noetherian regular scheme $X$ and a prime number $\ell$ invertible on $X$, the purity of Brauer groups holds for $X$ 
 with respect to the $\ell$-primary torsion part (\cite{RZ}, \cite{Th}, \cite{Fu}).
\item[(3)]
For a regular scheme $X$ over $\bF_p$, the purity of Brauer groups holds for $X$
 with respect to the $p$-primary torsion part (\cite{Mi2}, \cite{Gr1}, \cite{Sh}).
\item[(4)]
Let $p$ be a prime number, and let $k$ be a henselian discrete valuation field of characteristic zero whose residue field is perfect of characteristic $p$. Then for a smooth or semistable family $X$ over the integer ring of $k$, the purity of Brauer groups holds for $X$ with respect to the $p$-primary torsion part (\cite{Sat2} Corollary 4.5.2). This purity fact is a consequence of a result of Hagihara \cite{Sat2} Theorem A.2.6, which relies on the Bloch-Kato-Hyodo theorem on \'etale sheaves of $p$-adic vanishing cycles (\cite{bk}, \cite{Hy}).
\end{enumerate}
\end{rem}

\subsection{Milnor {\it K}-groups and filtration}\label{def:filtration}

\def\mK{\text{\it K}^{\text{\it M}}}
\def\mk{\text{\it k}^{\text{\it M}}}
\def\mh{\text{\it h}}
\def\ZZ{\text{\it Z}}

For a field $K$, we define the {\it Milnor K-group} $\mK_2(K)$ as
\[ \mK_2(K) :=(K^\times)^{\otimes 2}/J, \]
where $J$ denotes the subgroup of $(K^\times)^{\otimes 2}$ generated by symbols
\[ a_1 \otimes a_2 \quad \hbox{ with }\; a_1,a_2 \in K^\times \; \hbox{ and } a_1+a_2=0 \; \hbox{ or }1. \]
When $K$ is a discrete valuation field, we define the associated filtration $U^m\mK_2(K) \subset \mK_2(K)$ ($m \geq 0$) as the full group $\mK_2(K)$ if $m=0$, and as the subgroup generated by symbols of the form
\[ \{1+ \pi^m a, b\}\quad \hbox{ with } a \in \OO_K \hbox{ and } b \in K^\times \]
if $m \geq 1$, where $\OO_K$ denotes the integer ring of $K$ and $\pi$ denotes a prime element of $\OO_K$. Let $p$ be a prime number which is different from $\ch(K)$, and put
\[  H^2(K):=H^2(K,\mu_p^{\otimes 2}). \]
We define the filtration $U^m$ on $H^2(K)$ as that induced by $U^m$ on $\mK_2(K)$ via the norm residue symbol
\[ \mK_2(K) \lra H^2(K). \]
Now let $K$ be a henselian discrete valuation field of mixed characteristic $(0,p)$. Let $F$ be the residue field of $\OO_K$, and let $\pi$ be a prime element of $\OO_K$. Put
\[ \ZZ^q_F := \ker(d : \W^q_F \to \W^{q+1}_F), \qquad \W^q_{F,\log}:=\Image(\dlog : \mK_q(F) \to \W^q_F)). \]
Then by \cite{bk} Theorem (5.12), we have isomorphisms
\begin{equation}\label{gr:k2:1}
  \gr_U^mH^2(K)  \simeq
 \begin{cases}
 \W^2_{F,\log} \oplus \W^1_{F,\log}  & \hbox{($m=0$)},\\
 \W^1_F & \hbox{($0 < m < e'$, $p \hspace{-5pt}\not \vert m$)},\\
 \W^1_F/\ZZ_F^1 \oplus F/F^p  &\hbox{($0 < m < e'$, $p \vert m$)}
 \end{cases}
\end{equation}
defined by the following assignments respectively for $a_1,a_2 \in \OO_K$ and $b_1,b_2,b_3 \in \OO_K^\times$:
{\allowdisplaybreaks
\begin{align*} \{ b_1, b_2 \} + \{\pi, b_3 \} &  \mapsto (\dlog (\ol{b_1}) \wedge \dlog (\ol{b_2}), \dlog (\ol{b_3})), \\
  \{1+ \pi^ma_1, b_1 \} + \{1+ \pi^ma_2, \pi \} & \mapsto  \ol{a_1}  \cdot \dlog (\ol{b_1})-m^{-1} \cdot d\ol{a_2}\,, \\
  \{1+ \pi^ma_1, b_1 \} + \{1+ \pi^ma_2, \pi \} & \mapsto \big(\,\ol{a_1} \cdot \dlog (\ol{b_1}), \ol{a_2}\,\big).
\end{align*}}
\begin{lem}\label{lem:exact}
Let $K$ be a henselian discrete valuation field of mixed characteristic $(0,p)$. Let $F$ be the residue field of $\OO_K$. Then
\begin{enumerate}
\item[(1)]
There is a short exact sequence
\[ 0 \lra \Br(F) \oplus \H^1(F,\qz) \lra \Br(K) \lra \Br(K^\ur)^{G_F} \lra 0, \]
where $K^\ur$ denotes the maximal unramified extension of $K$ and $G_F=\Gal(\ol F/F)$ denotes the absolute Galois group of $F$.
\item[(2)]
If $K$ contains a primitive $p$-th root of unity, then the image of the map
\[ {}_p\Br(F) \otimes \mu_p(K) \lra {}_p\Br(K) \otimes \mu_p(K) \simeq H^2(K) \]
is contained in $U^{e'}H^2(K)$, where $e'$ denotes the natural number $p \cdot \ord_K(p)/(p-1)$.
\end{enumerate}
\end{lem}
\begin{pf}
(1) There is a short exact sequence
\[ 0  \lra \H^2(F,(K^\ur)^\times) \lra \H^2(K,\Gm) \lra \H^2(K^\ur,\Gm)^{G_F} \lra 0 \]
obtained from the Hochschild-Serre spectral sequence
\[ E_2^{u,v}=\H^u(F,\H^v(K^\ur,\Gm)) \Lra \H^{u+v}(K,\Gm), \]
and Hilbert's theorem 90. The assertion follows from the direct decomposition
\[ (K^\ur)^\times \simeq \bZ \times (\OO_K^\ur)^\times,\quad a \mapsto (\ord(a),a \pi^{-\ord(a)}), \]
where $\pi$ denotes a fixed prime element of $\OO_K$.
\par
(2) The assertion is a variant of \cite{bk} Lemma (5.1) (ii), whose details are left to the reader as an exercise.
\end{pf}

\section{Unramifiedness theorem for Brauer groups}\label{sect4}
\medskip
Let $k$ be a henselian discrete valuation field of characteristic $0$ whose residue field $\bF$ is finite and has characteristic $p$. Let $\OO_k$ be the integer ring of $k$ and put $S:=\Spec(\OO_k)$.

\subsection{Two generalizations on unramifiedness}\label{sect4-1}
Let $\cU$ be an integral scheme which is faithfully flat of finite type over $S$. Let $V$ be the divisor on $\cU$ defined by the radical of $(p)\subset \O_{\cU}$. We call $\w\in \Br(U)$ {\it unramified along} $V$, if $\w$ is contained in the image of $\Br(\cU) \ra \Br(U)$. If $\cU$ is regular, this condition is equivalent to that $\w$ belongs to the subgroup $\Br(\cU) \subset \Br(U)$. Following the ideas of Colliot-Th\'el\`ene--Saito in \cite{cts} \S2, we introduce two generalized notions of unramifiedness.
\begin{defn}\label{def:ur}
\begin{enumerate}
\item[(0)] We say that an \'etale morphism $f : B \to \cU$ is {\it quasi-cs along} $V$, if it satisfies the following two conditions.
\begin{enumerate}
\item[(i)] For any generic point $\eta$ of $V$, there exists exactly one connected component $B'$ of $B$ which splits completely over $\eta$ (namely, the image of $g:=f \vert_{B'}$ contains $\eta$ and $g^{-1}(\eta)$ is isomorphic to the sum of finitely many copies of $\eta$).
\item[(ii)] Each connected component of $B$ splits completely over some generic point of $V$.
\end{enumerate}
\item[(1)]
We say that $\w \in \Br(U)$ is {\it quasi-unramified along} $V$, if there exists an \'etale map $B \to \cU$ quasi-cs along $V$ such that $\w \vert_{B_k} \in \Br(B_k)$ belongs to the image of $\Br(B)$.
\item[(2)]
We say that $\w \in \Br(U)$ is $0$-{\it unramified}, if its specialization $\omega|_v \in \Br(v)$ is zero for any closed point $v$ on $U$ whose closure in $\cU$ is {\it finite over} $S$.
\end{enumerate}
\end{defn}
\begin{rem}\label{rem:ur}
\begin{enumerate}
\item[(1)]
Let $\omega \in \Br(U)$ be unramified along $V$. Then $\omega$ is quasi-unramified along $V$ obviously, and we see that $\omega$ is $0$-unramified as follows. Indeed, for a closed point $i: v \ra U$ whose closure in $\cU$ is finite over $S$, there is a commutative diagram of schemes
\[\xymatrix{ v \;\ar@<-1pt>@{^{(}->}[r] \ar[d]_i & \Spec(\OO_v) \ar[d] \\ U  \;\ar@<-1pt>@{^{(}->}[r] & \cU, }\]
where $\OO_v$ is the integer ring of $\kappa(v)$. Hence $\omega|_v$ is zero by the fact that $\Br(\OO_v)$ is zero.
\item[(2)]
For a generic point $\eta$ of $V$, let $A_{\eta}$ be the henselization of $\cO_{\cU\!,\eta}$ at its maximal ideal, and let $K_{\eta}$ be the fraction field of $A_{\eta}$. Then we have
\addtocounter{equation}{2}
\begin{equation}\label{isom:rem:1}
\bigoplus_{\eta \in V^0} \ \Br(A_{\eta}) \simeq \varinjlim_B \ \Br(B),
\end{equation}
where $B$ ranges over all \'etale $\cU$-schemes which are quasi-cs along $V$ (note that the set of such $B$'s endowed with a natural semi-order is co-filtered). Hence $\w \in \Br(U)$ is quasi-unramified along $V$ if and only if its restriction to $\bigoplus_{\eta \in V^0}~\Br(K_{\eta})$ belongs to the subgroup $\bigoplus_{\eta \in V^0}~\Br(A_{\eta})$.
\item[(3)]
If $\cU$ is regular and the purity of Brauer groups holds for $\cU$ (in the sense of Definition \ref{def:purity}), then $\omega \in \Br(U)$ is quasi-unramified along $V$ {\it if and only if} $\omega$ is unramified along $V$. Indeed, assuming the purity of Brauer groups, one can easily see that the restriction map
\[ \H^3_V(\cU,\Gm) \lra \bigoplus_{\eta \in V^0} \ \H^3_{\eta}(\cU_{\eta},\Gm) \]
 is injective, and that the restriction map
\[ \Br(U)/\Br(\cU) \lra \bigoplus_{\eta \in V^0} \ \Br(K_{\eta})/\Br(A_{\eta}) \]
 is injective as well.
\end{enumerate}
\end{rem}
\subsection{Unramifiedness theorem}\label{sect4-2}
Let $\cX$ be either an integral proper flat scheme over $S$ or its henselization at a closed point. Let $Y$ be the divisor on $\cX$ defined by the radical of $(p) \subset \cO_{\scX}$. Let $Z \subset \cX$ be a closed subscheme of pure codimension one with $Y \subset Z$. Let $Z_f$ be the union of the irreducible components of $Z$ which is flat over $S$, and put $\cU:=\cX-Z_f$. Put
\[ U:=\cX-Z,\qquad V:=Y \cap \cU, \qquad d:=\dim(\cX)=\dim(\cU). \]
By Remark \ref{rem:ur}\,(1), the following implications hold for elements of $\Br(U)$:
\[ \hbox{quasi-unramified along } V \Longleftarrow \hbox{unramified along } V \Lra 0\hbox{-unramified}. \]
The main result of this section is the following implication:
\begin{thm}\label{thm:key}
Assume that $\cU$ is regular in codimension $1$. If $\cX$ is henselian local, then assume further that $\cU=\cX$ {\rm(}and $V=Y=Z${\rm)} and that all irreducible components of $Y$ are regular. Then an arbitrary $0$-unramified element of $\Br(U)$ is quasi-unramified along $V$.
\end{thm}
\begin{rem}\label{rem:equiv}
\begin{enumerate}
\item[(1)]
By Remark \ref{rem:ur}\,(2), the assertion in Theorem \ref{thm:key} is equivalent to the claim that any $0$-unramified element maps to zero under the natural map \[ \Br(U) \lra \bigoplus_{\eta \in V^0} \ \Br(\Kh)/\Br(\Ah). \]
\item[(2)]
If $U$ is regular, then $\Br(U)$ is torsion. In this case the prime-to-$p$ part of Theorem \ref{thm:key} is due to Colliot-Th\'el\`ene--Saito \cite{cts} Th\'eor\`eme 2.1.
\end{enumerate}
\end{rem}
\noindent
By Remark \ref{rem:ur}\,(3), we obtain the following corollary:
\begin{cor}\label{cor:key}
Assume that the purity of Brauer groups holds for $\cU$. Then the following three conditions for $\omega \in \Br(U)$ are equivalent{\rm:}
\begin{enumerate}
\item[(1)] $\omega$ is $0$-unramified.
\item[(2)] $\omega$ is quasi-unramified along $V$.
\item[(3)] $\omega$ is unramified along $V$, i.e., belongs to $\Br(\cU)$.
\end{enumerate}
\end{cor}
\noindent
We will prove Theorem \ref{thm:key} in \S\S\ref{sect4-3}--\S\ref{sect4-5} below.

\subsection{An id\`ele class group}\label{sect4-3}
Let the notation be as in \S\ref{sect4-2}.  The following construction of id\`ele class groups is a slight modification of 
\cite{Ke} \S3 and \S4 (cf.\ \cite{ss:localring}).

\begin{defn}
\begin{enumerate}
\item[(1)]
A {\it chain} on $\cX$ is a sequence of points $P=(p_0,p_1,\dots,p_s)$ on $\cX$ such that 
\[ \overline{\{p_0\}}\subset\overline{\{p_1\}}\subset\cdots \subset \overline{\{p_s\}}, \]
where the closures are take on $\cX$. The dimension $d(P)$ of a chain $P=(p_0,\dots,p_s)$ is defined as $\dim\overline{\{p_s\}}$.
\item[(2)]
A {\it Parshin chain} on $(\cX,Z)$ is a chain $P=(p_0,\dots,p_s)$ such that $\dim \overline{\{p_i\}} =i$ for $0\leq i\leq s$ and such that $p_i \in Z$ for $i \leq s-1$ and $p_s\in U$.
\item[(3)]
A {\it {\it Q}-chain} on the pair $(\cX,Z)$ is a chain $P=(p_0,\dots,p_{s-2},p_s)$ such that $\dim \overline{\{p_i\}} =i$ for $i\in \{0,1,\dots,s-2,s\}$ and that $p_i\in Z$ for $i\leq s-2$ and $p_s\in U$.
\end{enumerate}
\end{defn}

For a chain $P=(p_0,\dots,p_s)$, let $\cO_{\sscX,P}^\h$ be the finite product of henselian local rings constructed as in \cite{KS} Definition 1.6.2\,(1) and \cite{Ke} Definition 3.1, and let $K_P$ be the product of its residue fields. By definition $K_P$ is a ring over the residue field of $p_s$ and we have the induced map if $p_s \in U$:
\begin{equation}\label{chain.eq}
\iota_P: \Spec(K_P) \lra U.
\end{equation}
Let $\cP$ (resp.\ $\cQ$) be the set of Parshin chains (resp.\ {\it Q}-chains) on  $(\cX,Z)$. For an integer $j>0$, let $\cP_j$ be  the set of $P\in \cP$ with $d(P)=j$). 

\begin{rem}\label{kerz.rem1}
A closed point $x\in U_0$ gives rise to $P_x=(y,x) \in \cP_1$, where $y$ is the unique point of $\overline{\{x\}}\cap Y$.
We identify $U_0$ with a subset of $\cP_1$ by the assignment $x \mapsto P_x$.
\end{rem}

\begin{defn}
For a Weil divisor $D$ such that $|D|\subset Z$, we define the {\it id\`ele group} of $(\cX\!,D)$ as
\[
I(\cX\!,D) := \Coker\bigg(\underset{P\in \cP_d}{\bigoplus}\ U^{D(P)}\mK_{d(P)-1}(K_P) \to 
\underset{P\in \cP}{\bigoplus}\ \mK_{d(P)-1}(K_P) \bigg),
\]
where $D(P)$ is the multiplicity of the prime
divisor $\overline{\{p_{d-1}\}}$ in $D$ for $P=(p_0,\dots,p_d)\in \cP_d$.
For Weil divisors $D,D'$ such that $|D|,|D'|\subset Z$ and $D'\geq D$,
we have a natural surjective map
\[
I(\cX\!,D') \to I(\cX\!,D).
\]
\end{defn}

There exists a natural homomorphism
\[ \psi: \underset{Q\in \cQ}{\bigoplus}\; \mK_{d(P)-1}(K_P)  \to I(\cX\!,D) \]
whose $(P,Q)$-componets $\psi_{P,Q}$ are defined as follows. Take $Q=(p_0,\dots,p_{s-2},p_s)\in \cQ$.
If $P=(p_0,\dots,p_{s-2},p_{s-1},p_s)\in \cP$ for $p_{s-1}\in Z$, 
the natural inclusion $K_Q\subset K_P$ induces 
(note $d(Q)=d(P)$)
\[
\psi_{P,Q}: \mK_{d(Q)-1}(K_Q)\to \mK_{d(P)-1}(K_P).
\]
If $P=(p_0,\dots,p_{s-2},p_{s-1})\in \cP$ for $p_{s-1}\in U$, 
$K_P$ is the product of the residue fields of discrete valuations on $K_Q$
induced by $p_{s-1}$. Thus we have the residue symbol
(note $d(Q)=d(P)+1$)
\[
\psi_{P,Q}: \mK_{d(Q)-1}(K_Q)\to \mK_{d(P)-1}(K_P).
\] 

\begin{defn}
Let $D$ be a Weil divisor such that $|D|\subset Z$. We define the {\it id\`ele class group} of $(\cX\!,D)$ as
\[
C(\cX\!,D)=\Coker\bigg(\underset{P \in \cQ}{\bigoplus}\; \mK_{d(P)-1}(K_P)
\rmapo{\psi}  I(\cX\!,D)\bigg).
\]
\end{defn}
\medbreak

Note that for $P\in \cP$, $K_P$ is a product of $d(P)$-dimensional local fields, and that there is a natural injective map due to Kato \cite{kk:cft} \S3.4 Proposition 3
\[ \Phi_P:\Br(K_P) \to  \Hom(\mK_{d(P)-1}(K_P),\qz). \]

\begin{prop}\label{recmap}
There exists a canonical homomorphism
\[
\Phi: \Br(U) \to \indlim {|D|\subset Z} \Hom(C(\cX\!,D),\qz)
\]
fitting into the following commutative diagram for any $P\in \cP${\rm:}
\[
\xymatrix{
 \Br(U) \ar[r]^-{\Phi} \ar[d]_{\iota_P^*} & \varinjlim_{|D|\subset Z} \ \Hom(C(\cX\!,D),\qz) \ar[d] \\
 \Br(K_P) \ar[r]^-{\Phi_P}& \Hom(\mK_{d(P)-1}(K_P),\qz)
}
\]
where the left vertical map is induced by the map \eqref{chain.eq} and 
the right vertical map is induced by
the natural map $\mK_{d(P)-1}(K_P)\to I(\cX\!,D)$.
\end{prop}
\begin{pf}
The assertion follows from the same argument as for the construction of the reciprocity map in \cite{KS} \S3 using the reciprocity law proved in loc.\ cit.\ (3.7.4).
\end{pf}

\begin{prop}\label{density}
For any Weil divisor $D$ such that $|D|\subset Z$, the natural map
\[ \bigoplus_{x\in U_0} \ \bZ \lra C(\cX\!,D) \]
is surjective {\rm(}cf. Remark {\rm\ref{kerz.rem1})}.
\end{prop}
\begin{pf}
The assertion is proved by the same argument as in \cite{Ke} Corollary 6.8 and the following fact (cf.\ \cite{Ke} Corollary 6.7). For an arbitrary $P\in \cP$, there exists an integer $m$ such that the natural map
\[ \mK_{d(P)-1}(K_P) \lra C(\cX\!,D) \]
annihilates $U^{m'} \mK_{d(P)-1}(K_P)$ for any $m' \ge m$.
\end{pf}

\subsection{Proof of Theorem \ref{thm:key}}\label{sect4-5}
Let the notation be as in \S\ref{sect4-2}. For $P=(p_0,\dots,p_d)\in \cP_d$, write $A_P=\cO_{\scX,P'}^\h$ with $P'=(p_0,\dots,p_{d-1})$ and let $F_P=K_{P'}$ be the product of the residue fields of $A_P$. For each $\eta \in Y^0$, let $\cR_\eta \subset \cP_d$ be the set of Parshin chains $P=(p_0,\dots,p_d)\in \cP_d$ with $p_{d-1}=\eta$. For $\eta\in Y^0$ and $P \in \cR_\eta$, there are natural injective ring homomorphisms
\[ \iota_{\eta,P}: A_\eta\to A_P,\quad K_\eta\to K_P.\]
We will prove the following lemma:
\begin{lem}\label{lem:inj}
For each $\eta \in V^0$, the maps $\iota_{\eta,P}$ with $P \in \cR_\eta$ induce an injective map
\stepcounter{equation}
\begin{equation}\label{res:zeta}
 \Br(K_\eta)/\Br(A_\eta) \lra \prod_{P\in \cR_\eta} \ \Br(K_P)/\Br(A_P).
\end{equation}
\end{lem}
We first prove the theorem admitting this lemma. Indeed, we get
\[ \Ker\bigg(\Br(U) \to \underset{x\in U_0}{\prod} \Br(x)\bigg) \subset \Ker\bigg(\Br(U) \to \underset{P\in \cP_d}{\prod} \Br(K_P)\bigg) \]
by Propositions \ref{density} and \ref{recmap} and the injectivity of $\Phi_P$. Hence Lemma \ref{lem:inj} implies the theorem.
\par\medskip
\def\Ash{A_{\ol {\eta}}}
\def\Ksh{K_{\ol {\eta}}}
\begin{pf*}{Proof of Lemma \ref{lem:inj}}
Fix an $\eta \in V^0$. Let $\Ash$ be the strict henselization of $\Ah$, and let $\Ksh$ be the fraction field of $\Ash$. We define $A_{\ol P}$ as the product of the strict henselizations of the direct factors of $A_P$ (i.e., a product of copies of $\Ash$), and define $K_{\ol P}$ as the product of the fraction fields of the direct factors of $A_{\ol P}$ (i.e., a product of copies of $\Ksh$). There is a commutative diagram with exact rows
\[ \xymatrix{ 0 \ar[r] &  \H^1(\eta,\qz) \ar[r] \ar[d]_a &  \Br(\Kh)/\Br(\Ah) \ar[r] \ar[d]_{\eqref{res:zeta}} &  \Br(\Ksh) \ar[d] \\ 0 \ar[r] & \displaystyle \prod_{P\in \cR_\eta} \ \H^1(F_P,\qz) \ar[r] & \displaystyle  \prod_{P\in \cR_\eta} \ \Br(K_P)/\Br(A_P) \ar[r] & \displaystyle  \prod_{P\in \cR_\eta}\ \Br(K_{\ol P}), }\]
where the vertical arrows are restriction maps, and the exactness of each row follows from Lemma \ref{lem:exact}\,(1). The right vertical map is obviously injective. It remain to show the injectivity of the left vertical map $a$, which we prove in what follows. Recall that $\k(\eta)$ is the function field of $Y_\eta=\overline{\{\eta\}} \subset \cX$, and that $Y_\eta$ is either a integral proper scheme over a finite field, or its henselization at a regular closed point. 

Take a dense open regular subset $Y' \subset Y_\eta$ if $Y_\eta$ is proper over a finite field, and let $Y' = Y_\eta$ if $Y_\eta$ is henselian local. Let $\cR'$ be the subset of $\cR_\eta$ consisting of all Parshin chains $P=(p_0,p_1,\dots,p_d)\in \cP_d$ with $p_i \in Y'$ for $i=0,1,\dotsc,d-1$. 
For $P=(p_0,p_1,\dots,p_d) \in \cR'$, let $Y_P$ be the henselization of $Y'$ along $(p_0,p_1,\dots,p_{d-2})$,
which is a direct sum of the spectra of discrete valuation rings.
Let $z_P$ be the direct sum of the closed points of $Y_P$ and put $x_P=p_0$.
There is a commutative diagram
\[\xymatrix{
0  \ar[r] & \H^1(Y',\qz) \ar[r] \ar[d]_c & \H^1(\eta,\Z/p\Z) \ar[r] \ar[d]_{a'} & \displaystyle \bigoplus_{y \in (Y')^1} \ \H^2_y(Y',\qz) \ar[d]_b \\
& \displaystyle \prod_{P \in \cR'} \ \H^1(x_P,\qz) \ar[r]^d & \displaystyle \prod_{P \in \cR'} \ \H^1(F_P,\qz) \ar[r] &\displaystyle \prod_{P \in \cR'} \ \H^2_{z_P}(Y_P,\qz),}\]
where the vertical arrows are natural restriction maps, and the upper row is exact by the purity of branch locus \cite{SGA2} Expos\'e X Th\'eor\`eme 3.4\,(i). The lower row is a complex (but not necessarily exact) and the arrow $d$ is defined as the composite map
\[ d_P : \H^1(x_P,\qz) \simeq \H^1(\Spec(\cO_{Y'\!,\hspace{1pt}p_0}^\h),\qz) \lra \H^1(Y_P,\qz) \to \H^1(F_P,\qz) \]
for each $P \in \cR'$. The map $d_P$ for $P=(p_0,p_1,\dots,p_d) \in \cR'$ is injective if the closures $\overline{\{p_i\}} \subset Y'$ are regular at $x_P=p_0$ for $0 \le i \le d-2$. Therefore, in order to show that $a'$ is injective, it is enough to verify that $b$ and $c$ are injective, by a simple diagram chase in the above diagram. The map $c$ is injective by the \v{C}ebotarev density theorem \cite{Se2} Theorem 7, if $Y_\eta$ is proper over a finite field. If $Y'$ is henselian local, then $c$ is obviously injective. The injectivity of $b$ is checked as follows. We have
\[ \qquad\; \H^2_y(Y',\qz) \simeq  \H^0(y,R^2i_y^!\qz) \qquad \hbox{($i_y:y \hra \Spec(\cO_{Y,y}^\h)$)}\]
for any $y \in (Y')^1$, and we have
\[ \H^2_{z_P}(Y_P,\qz) \simeq  \H^0(z_P,R^2i_P^!\qz)  \qquad \hbox{($i_P : z_P \hra Y_P$)}\]
for any $P \in \cR'$. Moreover if $P=(p_1,p_2\dotsc,p_d)$ with $p_{d-1}=y$, then we have $R^2i_P^!\qz \simeq \tau^*R^2i_y^!\qz$, where $\tau:z_P\to y$ is the natural map which is essentially \'etale.
Hence $b$ is injective. Thus $a'$ and $a$ are injective and we obtain Lemma \ref{lem:inj}, which completes the proof of Theorem \ref{thm:key}.
\end{pf*}

\subsection{Application to arithmetic schemes}\label{sect4-6}
Let $\cX$ be an integral scheme which is {\it proper} flat of finite type over $\Spec(\Z)$.
Put $\gfiber:=\cX \otimes_{\Z}\Q$.
\begin{thm}
Assume that the purity of Brauer groups holds true for $\cX$. Then there is an exact sequence{\rm:}
\[ 0 \lra \Br(\cX) \lra \Br(\gfiber) \lra \prod_{x \in \gfiber_0} \ \Br(x)/\Br(\fo_x), \]
where for $x \in \gfiber_0$, $\fo_x$ denotes the integral closure of $\Z$ in $\kappa(x)$.
\end{thm}
\begin{rem}
For $x \in \gfiber_0$, $\Br(\fo_x)$ is a finite $2$-torsion group by the classical Hasse principle.
\end{rem}
\begin{pf}
The above sequence is a complex by the properness of $\cX$. We show that the resulting specialization map
\[ \Br(\gfiber)/\Br(\cX) \lra \prod_{x \in \gfiber_0} \ \Br(x)/\Br(\fo_x) \] is injective.
Let $P$ be the set of all prime numbers. For $p \in P$, let $\Z_p^\h$ be the henselization of $\Z$ at $(p)$. Put $\Q_p^\h:=\Frac(\Z_p^\h)$, and let $Q_p$ be the set of all closed points on $\gfiber_{\Q_p^\h}$. We construct the map
\[ \alpha:\prod_{x \in X_0} \ \Br(x)/\Br(\fo_x) \lra \prod_{p \in P} \ \prod_{v \in Q_p} \ \Br(v) \]
as follows. Let $x$ be a closed point on $\gfiber$ and let $v$ be a closed point on $\gfiber_{\Q_p^\h}$.
We define the $(x,v)$-component of $\alpha$ as the natural restriction map (resp.\ the zero map), if the composite map $v \ra \gfiber_{\Q_p^\h} \ra \gfiber$ factors through $x \ra \gfiber$ (resp.\ otherwise).
Note that for $v \in Q_p$, there exists a unique $x \in X_0$ such that the composite map $v \ra \gfiber_{\Q_p^\h} \ra \gfiber$ factors through $x \ra \gfiber$, and that this uniqueness implies the well-definedness of $\alpha$.
Now let us consider a commutative diagram of specialization maps
\[\xymatrix{ \Br(\gfiber)/\Br(\cX) \ar[r] \ar[d] & \displaystyle  \prod_{x \in X_0} \ \Br(x)/\Br(\fo_x) \ar[d]^\alpha \\
\displaystyle  \bigoplus_{p \in P} \ \Br\big(\gfiber_{\Q_p^\h}\big)\big/\Br\big(\cX_{\Z_p^\h}\big) \ar[r] & \displaystyle  \prod_{p \in P} \ \prod_{v \in Q_p}\ \Br(v), }\] whose commutativity follows from the definition of $\alpha$. In this diagram, the bottom horizontal arrow is injective by Corollary \ref{cor:key}, and the left vertical arrow is injective by the localization exact sequences
{\allowdisplaybreaks
\begin{align*}\xymatrix{\dotsb \ar[r] & \Br(\cX) \ar[r] & \Br(\gfiber) \ar[r] & \displaystyle \bigoplus_{p \in P} \ \H^3_{Y_p}(\cX,\Gm) \ar[r] & \dotsb,} \\ \xymatrix{ \dotsb \ar[r] & \Br\big(\cX_{\Z_p^\h}\big) \ar[r] & \Br\big(\gfiber_{\Q_p^\h}\big) \ar[r] & \H^3_{Y_p}\big(\cX_{\Z_p^\h},\Gm\big) \ar[r] & \dotsb},
\end{align*}
}with $Y_p:=\cX \otimes_{\Z}\bF_p$, and the excision isomorphism \[ \H^3_{Y_p}(\cX,\Gm) \simeq \H^3_{Y_p}(\cX_{\Z_p^\h},\Gm) \quad \hbox{ for each \;$p \in P$.} \] Hence the upper horizontal map is injective as well and we obtain the theorem.
\end{pf}

\section{zero-cycles on cubic surfaces}\label{sect.cubic}
\medskip
In this section, we compute $\A_0$ of cubic surfaces explicitly using the unramifiedness theorem proved in the previous section.
\subsection{Setting and results}
Let $k,\OO_k,\bF$ and $p$ be as in the beginning of \S\ref{sect4}.
Let $a$ be an element of $k^\times$ which is not a cube in $k$.
We are concerned with a cubic surface
\[ X := \{ T_0^3+T_1^3+T_2^3+ a T_3^3=0 \} \subset  \Proj(k[T_0,T_1,T_2,T_3]) = \bP^3_k. \]
Let $\zeta_3$ be a primitive cubic root of unity in $\ol k$.
\begin{thm}\label{thm.cubic}
\begin{itemize}
\item[(1)]
Assume $p\not=3$. Then we have
\[
\A_0(X)\simeq
\left.\left\{\begin{gathered}
 0 \\
 \bZ/3 \\
 \bZ/3 \oplus \bZ/3 \\ 
\end{gathered}\right.\quad
\begin{aligned}
&\text{if $\ord_k(a)\equiv 0\mod 3$,}\\
&\text{if $\ord_k(a)\not\equiv 0\mod 3$ and $\zeta_3 \not\in k$,} \\
&\text{if $\ord_k(a)\not\equiv 0\mod 3$ and $\zeta_3 \in k$.}
\end{aligned}\right.
\]
\item[(2)]
Assume $p=3$, $\ord_k(a) \equiv 1\mod 3$ and $\zeta_3 \in k$. Then we have
\[ \A_0(X) \simeq \bZ/3 \oplus \bZ/3. \]
\end{itemize}
\end{thm}
(1) is stated in \cite{cts} Example 2.8 under a slightly simpler setting. We include a proof of (1) here for the convenience of the reader. (2) is a new result and would be the first example of a potentially rational surface which splits over a wildly ramified extension and whose $\A_0$ is computed explicitly. It would be interesting to find cycles which generate $\A_0(X)$ in the theorem. \par To prove Theorem \ref{thm.cubic}, we need the following three facts, where $X$ is as before.
\begin{prop}[{{\bf Colliot-Th\'el\`ene}}]\label{prop.ex1}
The following map induced by the Brauer-Manin pairing is injective{\rm:}
\[ \A_0(X) \lra \Hom(\Br(X)/\Br(k),\qz). \]
\end{prop}
\begin{pf}
The case $X(k) \ne \emptyset$ is stated in \cite{CT1} Proposition 5. Otherwise, the assertion follows from his injectivity result in loc.\ cit.\ Proposition 7\,(b) and the same arguments as in loc.\ cit.\ Proposition 5 (cf.\ \cite{Bl} Theorem (2.1), Proposition (A.1)).
See also \cite{ss:cycle} Theorem A and \cite{kahn} p.\ 70 Corollaire 2 for generalizations.
\end{pf}
\begin{prop}\label{prop.ex2}
Let $\cU$ be a regular scheme which is faithfully flat over $S$ and satisfies $\cU \otimes_{\OO_k} k \simeq X$.
Let $\eta$ be a generic point of $\cU \otimes_{\OO_k} \bF$, let $\Ah$ be the henselization of $\cO_{\hspace{-1pt}\cU\hspace{-1.3pt},\eta}$ and let $\Kh$ be the fraction field of $\Ah$.
Assume that
\addtocounter{equation}{3}
\begin{equation}\label{eq1.prop.ex2}
 \iota^{-1}(\Br(\Ah)) \subset \Image(\Br(k) \to \Br(X)), \quad \hbox{ where } \; \iota : \Br(X) \to \Br(\Kh). \end{equation}
Then the following map induced by the Brauer-Manin pairing is surjective{\rm:}
\begin{align*} \A_0(X) & \lra \Hom(\Br(X)/\Br(k),\qz),\end{align*}
\end{prop}
\begin{pf}
By \eqref{eq1.prop.ex2}, Theorem \ref{thm:key} and Remark \ref{rem:ur}\,(2), the map in question has dense image.
As we mentioned in \S\ref{sect1-1'}, $\Br(X_L)/\Br(L)$ is zero for $L=k(\sqrt[3]{a})$.
Hence $\Br(X)/\Br(k)$ is a finite $3$-torsion and the assertion follows.
\end{pf}
To state the third fact, we assume that $k$ contains $\zeta_3$, and fix an isomorphism
\begin{equation}\label{eq1.thm.cubic}
{}_3\Br(k(X)) \simeq H^2(k(X),\mu_3^{\otimes 2}), \quad x \mapsto x \otimes \zeta_3,
\end{equation}
where $k(X)$ is the function field of $X$. Consider rational functions \[ f :=\frac{T_0+\zeta_3 T_1}{T_0+T_1},\;\; g :=\frac{T_0+T_2}{T_0+T_1} \in k(X)^\times \] and put \[ \bs{e}_1=(a,f)_{\zeta_3},\; \bs{e}_2=(a,g)_{\zeta_3} \in {}_3\Br(k(X)), \]
where for $u,v \in k(X)^\times$, $(u,v)_{\zeta_3}$ denotes the inverse image of $\{u,v\} \in H^2(k(X),\mu_3^{\otimes 2})$ under the isomorphism \eqref{eq1.thm.cubic}.
\addtocounter{thm}{2}
\begin{thm}[{{\bf Manin \cite{M2}}}]\label{prop.manin}
$\bs{e}_1$ and $\bs{e}_2$ belong to $\Br(X)$, and $\Br(X)/\Br(k)$ is a free $\bZ/3$-module of rank $2$ generated by $\bs{e}_1$ and $\bs{e}_2$.
\end{thm}
\def\cOXeta{\cO_{\hspace{-1.5pt}\cX\hspace{-1.5pt},\eta}}
\subsection{Proof of Theorem \ref{thm.cubic} (1)}
Without loss of generality, we may assume that $a \in \OO_k-\{ 0 \}$ with $\ord_k(a)=0,1$ or $2$. Consider a projective flat model of $X$ over $S:=\Spec(\OO_k)$
\begin{equation}\label{eq0.thm.cubic}
\cX=\Proj\big(\OO_k[T_0,T_1,T_2,T_3]/(T_0^3+T_1^3+T_2^3+ a T_3^3)\big) \subset \bP^3_S.
\end{equation}
Let $\eta$ be the generic point of $Y:=\cX \otimes_{\OO_k} \bF$, let $\Ah$ be the henselization of $\cOXeta$ and let $\Kh$ be the fraction field of $\Ah$.
We divide the problem into 4 cases as follows. \par
\medskip
\noindent
{\bf Case (i):} $p\not=3$, $a\in \OO_k^\times$ and $\zeta_3 \in k$.
In this case it is easy to see that $\cX$ is smooth over $S$. Once we show  \begin{equation}\label{eq1-01.thm.cubic}\Br(X)=\Br(\cX)+\Image(\Br(k)\to \Br(X)), \end{equation}
 then $\A_0(X)=0$ by Proposition \ref{prop.ex1}. Since we have
\[ \Br(\cX) = \ker (\Br(X) \to \Br(K_\eta)/\Br(A_\eta)), \]
by the purity of Brauer groups for $\cX$, it is enough to show \begin{equation}\label{eq1-00.thm.cubic} \Image(\Br(X)\to \Br(K_\eta))\subset \Br(A_\eta)+\Image(\Br(k)\to \Br(K_\eta)), \end{equation}
in order to show \eqref{eq1-01.thm.cubic}. Since $p \ne 3$ and $\zeta_3 \in k$, there is an exact sequence
\begin{equation}\label{eq1-0.thm.cubic} 0\lra {}_3\Br(A_\eta) \otimes \mu_3(k) \lra H^2(k(X),\mu_3^{\otimes 2}) \rmapo{\delta} \k(\eta)^\times/3 \lra 0 \end{equation}
and it is easy to check $\delta(\bs{e}_i\otimes \zeta_3)=0$ for $i=1,2$, which implies \eqref{eq1-00.thm.cubic} by Theorem \ref{prop.manin}.
\par\medskip
\noindent
{\bf Case (ii):} $p\not=3$, $a\in \OO_k^\times$ and $\zeta_3 \not\in k$.
In this case, the assertion is reduced to the case (i) immediately by a standard norm argument.
\par\medskip
\noindent
{\bf Case (iii):} $p\not=3$, $\ord_k(a)=1,2$ and $\zeta_3 \in k$.
One can easily check that the fixed model $\cX$ is regular at $\eta \in Y$, i.e.,
$\cOXeta$ is a discrete valuation ring.
By Propositions \ref{prop.ex1}, \ref{prop.ex2} and Theorem \ref{prop.manin}, we have only to show
\[ \iota^{-1}(\Br(\Ah)) \subset \Image(\Br(k) \to \Br(X)), \quad \hbox{ with } \; \iota : \Br(X) \to \Br(\Kh). \]
Note that we have the exact sequence \eqref{eq1-0.thm.cubic} in this case as well. For \[ \bs{\omega}=\alpha \, \bs{e}_1+ \beta \, \bs{e}_2\in {}_3\Br(X)\subset {}_3\Br(k(X)) \quad (\alpha,\beta\in \bZ/3), \] we have \[ \delta(\bs{\omega} \otimes \zeta_3)=\ord_k(a) \cdot f^\alpha g^\beta = \ord_k(a) \cdot \bigg(\frac{T_0+\zeta_3 T_1}{T_0+T_1}\bigg)^\alpha\bigg(\frac{T_0+T_2}{T_0+T_1}\bigg)^\beta  \in \k(\eta)^\times/3, \] where we regarded $f$ and $g$ as rational functions on $Y$. Thus it is enough to show
\addtocounter{thm}{4}
\begin{lem}\label{claim2.thm.cubic}
For $\alpha,\beta \in \bZ$, assume that $f^\alpha g^\beta \in \k(\eta)^\times$ belongs to $(\k(\eta)^\times)^3$. Then we have $\alpha \equiv \beta \equiv 0 \mod 3$.
\end{lem}
\begin{pf}
Let $E$ be the elliptic curve over $\bF$ defined as
\[ E:=\Proj\big(\bF[T_0,T_1,T_2]/(T_0^3+T_1^3+T_2^3)\big)\subset \bP^2_\bF. \]
It is easy to see that $\k(\eta)$ is the rational function field in one variable over the function field $\bF(E)$.
Since $(\bF(Y)^\times)^3\cap \bF(E)^\times =(\bF(E)^\times)^3$ and $f,g \in \bF(E)^\times$, 
the assumption of the lemma implies that $f^\alpha g^\beta \in (\bF(E)^\times)^3$. We now look at the divisors on $E$
\[ \div_E(f)=3([P]-[O]), \qquad \div_E(g)=3([Q]-[O]), \]
where we put
\[ O:=(1:-1:0),\; P:=(1:-\zeta_3:0),\; Q:=(1:0:-1). \]
Take $O$ to be the origin of the elliptic curve $E$, and define zero-cycles $C,C'$ on $E$ as $C:=[P]-[O]$, $C':=[Q]-[O]$.
Since $f^\alpha g^\beta = h^3$ for some $h \in \k(\eta)^\times$ by assumption, we have
\[ \alpha\cdot C +\beta\cdot C' = \div_E(h)\]
as zero-cycles, and the residue class $\alpha \cdot\ol{C}+ \beta \cdot\ol{C'}$ is zero in $\A_0(E)\simeq E(\bF)$. Hence the assertion follows from the linear independence of $\ol{C}$ and $\ol{C'}$ in the $\bZ/3$-vector space ${}_3E(\bF)$. \end{pf}
\par\medskip
\noindent
{\bf Case (iv):} $p\not=3$, $\ord_k(a)=1,2$ and $\zeta_3 \not\in k$.
Consider the scalar extension
\[ \cX \otimes_{\OO_k} \OO_L \lra \cX \qquad \hbox{($L:=k(\zeta_3)$)}, \]
which is \'etale by the assumption $p \ne 3$. Then we have \[ \iota^{-1}(\Br(\Ah)) \subset \Image(\Br(k) \to \Br(X)), \quad \hbox{ with } \; \iota : \Br(X) \to \Br(\Kh) \] by the previous case and a standard norm argument. Therefore by Propositions \ref{prop.ex1} and \ref{prop.ex2}, the assertion is reduced to the following proposition due to Colliot-Th\'el\`ene \cite{ct:letter}, which holds without the assumption $p \ne 3$:
\begin{prop}[{{\bf Colliot-Th\'el\`ene}}]\label{prop:CT}
Let $\bs{e}_i$\,{\rm($i=1,2$)} be the elements of $\Br(X_L)$ in Proposition {\rm\ref{prop.manin}}. Then
$\Br(X)/\Br(k)$ is a free $\bZ/3$-module of rank $1$ generated by $\Cor_{X_L/X} (\bs{e}_1)$, where $\Cor_{X_L/X}$ denotes the corestriction map $\Br(X_L) \to \Br(X)$.
\end{prop}
\begin{pf}
Put $G:=\Gal(L/k)$, which has order $2$. Let $\sigma$ be the generator of $G$. We prove
\addtocounter{equation}{2}
\begin{equation}\label{eq:CT}
\sigma(\bs{e}_1) = \bs{e}_1 \quad \hbox{ and } \quad \sigma(\bs{e}_2) = -\bs{e}_2 \quad \hbox{ in } \; {}_3\Br(X_L),
\end{equation}
which implies the assertion by a standard norm argument.
To prove \eqref{eq:CT}, we work with Galois cohomology groups of the function field \[ F:=L(X_L)=k(X)(\zeta_3). \]
Since ${}_3\Br(X_L) \subset H^2(F,\mu_3)$,
 it is enough to show the following two claims:
\begin{enumerate}
\item[(1)]
{\it For $u,v \in F^\times$, we have $\sigma((u,v)_{\zeta_3}) = -(\sigma(u),\sigma(v))_{\zeta_3}$ in $H^2(F,\mu_3)$.}
\item[(2)]
{\it We have} \[ \bigg\{ a, \frac{T_0+\zeta_3^{-1}T_1}{T_0 +T_1} \bigg\}=-\bigg\{a,\frac{T_0+\zeta_3 T_1}{T_0+T_1}\bigg\} \quad \hbox{\it in } \;\; H^2(F,\mu_3^{\otimes 2}). \]
\end{enumerate}
We first show (1). Since $\bZ/3 \simeq \mu_3^{\otimes 2}$ as $G_k$-modules, we have
\begin{align*} \sigma((u,v)_{\zeta_3}) & = \sigma(\{u,v\} \otimes \zeta_3)= \{\sigma(u), \sigma(v)\} \otimes  \sigma(\zeta_3) \\
 & = \{\sigma(u), \sigma(v)\} \otimes  \zeta_3^{-1} = - (\sigma(u),\sigma(v))_{\zeta_3}. \end{align*}
We next show (2). Take an affine open subset of $X$ as follows:
\[ \{ x^3+y^3+z^3+a = 0 \} \subset \bA^3_k \quad 
\bigg(x=\frac{T_0}{T_3},\; y=\frac{T_1}{T_3},\; z=\frac{T_2}{T_3}\bigg). \]
Then noting that $H^2(F,\mu_3^{\otimes 2})$ is a $3$-torsion, we compute
\begin{align*} \bigg\{ a, \frac{x+\zeta_3^{-1}y}{x +y} \bigg\}+\bigg\{a,\frac{x+\zeta_3 y}{x+y}\bigg\}
 & = \bigg\{ a, \frac{x^2-xy+y^2}{(x +y)^2} \bigg\} = \bigg\{ a, \frac{x^3+y^3}{(x +y)^3} \bigg\} \\
& = \{ a, -a-z^3\} = \bigg\{\frac{a}{\,z^3\,}, -\frac{a}{\,z^3\,}-1\bigg\} = 0.
\end{align*} 
This completes the proof of Proposition \ref{prop:CT} and Theorem \ref{thm.cubic}\,(1).
\end{pf}
\medskip
\def\ep{\epsilon}
\subsection{Proof of Theorem \ref{thm.cubic}\,(2)}
Without loss of generality, we may assume $a=\pi$ (a prime of $\OO_k$).
Let $\cX$ be the projective flat model of $X$ over $\OO_k$ defined in \eqref{eq0.thm.cubic}. We will use the following affine open subset:
\begin{equation}\label{eq0-1.thm.cubic}
U := \Spec(\OO_k[x,y,z]/(x^3+y^3+z^3+\pi)) \quad \bigg(x=\frac{T_0}{T_3},\; y=\frac{T_1}{T_3},\; z=\frac{T_2}{T_3}\bigg).
\end{equation}
Let $\cX_s$ be the special fiber of $\cX \to S$, which is irreducible. Let $Y$ be the reduced part of $\cX_s$ and let $\eta$ be the generic point $Y$. It is easy to see that we have
\begin{equation}\label{eq2.thm.cubic}
U \cap Y=\Spec(\bF[x,y,z]/(x+y+z)).
\end{equation}
Put $e:=\ord_k(3)$. We show the following lemma.
\addtocounter{thm}{2}
\begin{lem}\label{lem2.thm.cubic}
$\cOXeta$ is a discrete valuation ring with absolute ramification index $3e$. 
\end{lem}
\begin{pf}
Put $t=x+y+z$. Then we have
\stepcounter{equation}
\begin{equation}\label{eq2.0.thm.cubic}
 t^3+\pi(1+\pi^{e-1}u) = x^3+y^3+z^3+\pi = 0 \quad \hbox{ in } \;\; \cO_{U,\eta}=\cOXeta,
\end{equation}
where we put
\[ u :=-\ep_1 \big(x^2(y+z)+y^2(z+x)+z^2(x+y)+2xyz\big) \;\; \text{ with } \; \ep_1 :=3\pi^{-e} \in \OO_k^\times. \]
In view of \eqref{eq2.thm.cubic}, this implies that $t$ generates the maximal ideal of $\cOXeta$ (and that $u$ belongs to $\cOXeta^\times$).  Hence $\cOXeta$ is a noetherian one-dimensional local ring with maximal ideal generated by $t$, which is a discrete valuation ring. The ramification index of $\cOXeta$ over $\OO_k$ is $3$ by \eqref{eq2.0.thm.cubic}, which implies that the absolute ramification index of $\cOXeta$ is $3e$.
\end{pf}
Let $\Ah$ be the henselization of $\cOXeta$ as before, and let $K_\eta$ be its fraction field. We put
\[ H:=H^2(K_\eta,\mu_3^{\otimes 2}) \]
in what follows. Since $k$ contains $\zeta_3$, we have \[ {}_3\Br(\Ah) \otimes \mu_3(k) \subset U^{3e'}H \] by Lemma \ref{lem:exact}\,(2), where $U^*H$ denotes the filtration on $H$ defined in \S\ref{def:filtration} and $e'$ denotes $pe/(p-1)=3e/2$.
As before, by Propositions \ref{prop.ex1}, \ref{prop.ex2} and Theorem \ref{prop.manin}, we ought to show
\[ \iota^{-1}(\Br(\Ah)) \subset \Image(\Br(k) \to \Br(X)), \quad \hbox{ with } \; \iota : \Br(X) \to \Br(\Kh). \]
It is enough to show the following:
\addtocounter{thm}{1}
\begin{prop}\label{claim3.thm.cubic} Let $\alpha,\beta \in \bZ/3$, and put
\[ \bs{\omega} := \alpha \bs{e}_1 + \beta \bs{e}_2 = \alpha (\pi,f)_{\zeta_3} + \beta(\pi,g)_{\zeta_3} \in {}_3 \Br(X). \] Assume that
\[ \iota(\bs{\omega})\otimes\zeta_3=\alpha\{\pi,f\}+\beta\{\pi,g\}\in {}_3 \Br(K_\eta)\otimes \mu_3(k) \simeq H \]
 belongs to $U^{3e'}H$. Then we have $\alpha=\beta=0$, i.e., $\bs{\omega}=0$.
\end{prop}
Note that $\k(\eta)=\bF(y,z)$ by \eqref{eq2.thm.cubic}.
By \eqref{eq2.0.thm.cubic}, we have $\pi = -t^3(1+\pi^{e-1}u)^{-1}$ and
\[ \bs{\omega} = -\alpha \{1+\pi^{e-1}u,f\} - \beta \{1+\pi^{e-1}u, g\}  \in U^{3(e-1)}H. \]
One can derive the proposition easily from the following lemma, where the Bloch-Kato isomorphisms (see \eqref{gr:k2:1}) are defined with respect to the prime element $t \in \Ah$.
\begin{lem}\label{claim3-0.thm.cubic}
\begin{enumerate}
\item[(1)] $\{1+\pi^{e-1}u,f\}$ belongs to $U^{3e'-2}H$, whose residue class in $\gr_U^{3e'-2}H \simeq \Omega^1_\eta$ is $zdy-ydz$ up to the multiplication by a constant in $\bF^\times$.
\item[(2)] $\{1+\pi^{e-1}u,g\}$ belongs to $U^{3e-2}H$, whose residue class in $\gr_U^{3e-2}H \simeq \Omega^1_\eta$ is $zdy-ydz$ up to the multiplication by a constant in $\bF^\times$.
\end{enumerate}
\end{lem}
\def\te{{\tilde{e}}}
\begin{pf*}{\it Proof of Lemma \ref{claim3-0.thm.cubic}}
We first prove (1).
Since $x^3+y^3+z^3+\pi=0$, we have
\addtocounter{equation}{2}
\begin{equation}\label{eq3.1.thm.cubic}
u= \ep_1(2y^2z+2yz^2+tv) \quad \hbox{with}  \;\; v:=-(x^2+(t-x)^2+\pi \cdot t^{-1})
\end{equation}
(note that $\pi \cdot t^{-1}$ is contained in the maximal ideal of $\Ah$). On the other hand, we have
\begin{equation}\label{eq3.2.thm.cubic}
f=\frac{\,x+\zeta_3 y\,}{x+y}=1+(\zeta_3-1)\frac{y}{\,x+y\,} = 1+t^{e'}\frac{\,\ep_2 y\,}{z} + t^{e'+1}\frac{\ep_2 y}{z(z-t)},
\end{equation}
where $\ep_2:=(1-\zeta_3)\cdot t^{-e'} \in \Ah^\times$.
We define a filtration $U^n\Kh^\times$\,($n \ge 0$) on $\Kh^\times$ as the full group $\Kh^\times$ for $n = 0$ and the subgroup $\{1+t^nc\,| \; c \in \Ah \}$ for $n \ge 1$. We recall the following standard facts:
\addtocounter{thm}{2}
\begin{sublem}\label{sublem1.thm.cubic}
For $a,b \in \Ah$, $h \in U^n\Kh^\times$\,{\rm ($n \ge 0$)} and integers $\ell,m,\nu \ge 1$, we have
\begin{align*}
\tag{1}
 &\{1+t^\ell a+t^m b, h \} \equiv \{1+t^\ell a, h \} + \{1+t^m b, h \} \mod U^{\ell+m+n}H \\
\tag{2}
 &\{1+t^\ell a(1+t^m b)^{\pm \nu}, h\} \equiv \{1+t^\ell a, h\} \mod U^{\ell+m+n}H \\
\tag{3}
 &\{1+t^\ell a,1+t^m b \} \equiv -\bigg\{1+\frac{t^{\ell+m}ab}{\,1+t^\ell a\,},-t^\ell a \bigg\} \mod U^{\ell+2m}H.
\end{align*}
\end{sublem}
\begin{pf*}{\it Proof of Sublemma \ref{sublem1.thm.cubic}} (1) follows from the equality
\[ \frac{(1+t^\ell a)(1+t^m b)}{1+t^\ell a+t^m b}= 1+\frac{t^{\ell + m}ab}{1+t^\ell a+t^m b} \]
and (3) (cf.\ \cite{bk} Lemma (4.1)). The assertion (2) follows from a similar computation.
(3) follows from similar computations as in \cite{Sat2} Lemma 8.7.4. \end{pf*}
\noindent
We return to the proof of Lemma \ref{claim3-0.thm.cubic}\,(1), and put $\ep_3 := 3t^{-3e} \in \Ah^\times$. Then we have
{\allowdisplaybreaks \begin{align*}
&\{1+\pi^{e-1}u,f\}
= \bigg\{ 1 + \pi^{e-1}\ep_1(2y^2z+2yz^2+tv),1+t^{e'}\frac{\,\ep_2 y\,}{z} + t^{e'+1}\frac{\ep_2 y}{z(z-t)}\bigg\} \\
& \quad\equiv \bigg\{ 1 + t^{3e-2}\ep_3v, 1+t^{e'}\frac{\,\ep_2 y\,}{z} \bigg\}
 + \bigg\{ 1 + 2t^{3e-3}\ep_3(y^2z+yz^2),1+ t^{e'+1}\frac{\ep_2 y}{z(z-t)}\bigg\} \\
&\quad \phantom{==} + \bigg\{ 1 + 2t^{3e-3}\ep_3 y^2z,1+t^{e'}\frac{\,\ep_2 y\,}{z} \bigg\}
 + \bigg\{ 1 + 2t^{3e-3}\ep_3 yz^2,1+t^{e'}\frac{\,\ep_2 y\,}{z} \bigg\},
\end{align*}
}where the congruity holds modulo $U^{3e'-1}H$ by Sublemma \ref{sublem1.thm.cubic}\,(1)--(3) and \eqref{eq2.0.thm.cubic}. Therefore, it is enough to show
\begin{lem}\label{lem.symbol}
\begin{enumerate}
\item[(1)] The sum
\[ \bigg\{ 1 + t^{3e-2}\ep_3v, 1+t^{e'}\frac{\,\ep_2 y\,}{z} \bigg\}
 + \bigg\{ 1 + 2t^{3e-3}\ep_3(y^2z+yz^2),1+ t^{e'+1}\frac{\ep_2 y}{z(z-t)}\bigg\} \in U^{3e'-2}H \] maps to $zdy-ydz$ up to the multiplication by a constant in $\bF^\times$, under the Bloch-Kato isomorphism $\gr_U^{3e'-2}H \simeq \Omega^1_\eta$ with respect to the prime element $t \in \Ah$.
\smallskip
\item[(2)]
$\displaystyle \bigg\{ 1 + 2t^{3e-3}\ep_3y^2z,1+t^{e'}\frac{\,\ep_2 y\,}{z} \bigg\}$ belongs to $U^{3e'-1}H$.
\smallskip
\item[(3)]
$\displaystyle \bigg\{ 1 + 2t^{3e-3}\ep_3yz^2,1+t^{e'}\frac{\,\ep_2 y\,}{z} \bigg\}$ belongs to $U^{3e'}H$.
\end{enumerate}
\end{lem}
\begin{pf}
(1) Let $\ol{\ep_2\ep_3} \in \kappa(\eta)^\times$ be the residue class of $\ep_2\ep_3=3(1-\zeta_3)t^{-3e'} \in \Ah^\times$, which in fact belongs to $\bF^\times$ by \eqref{eq2.0.thm.cubic}. We note that $3|e'$ and that $\ch(\eta)=3$. By Sublemma \ref{sublem1.thm.cubic}\,(3) and \eqref{eq3.1.thm.cubic}, the first term (resp.\ the second term) of the sum of symbols in question maps to
\[ \ol{\ep_2\ep_3}\cdot\frac{(y+z)^2(zdy - ydz)}{z^2}
  \quad \hbox{$\bigg($resp.\ \; $\ol{\ep_2\ep_3}\cdot \dfrac{(2y^2z+yz^2)dy+(2y^2z+y^3)dz}{z^2}\bigg)$}\] 
under the Bloch-Kato isomorphism, cf.\ \eqref{gr:k2:1}. Since $\ch(\eta)=3$, the sum of these 1-forms agrees with $\ol{\ep_2\ep_3}(zdy-ydz)$.
\par
(2) Noting that $e' \ge 3$ and that $H$ is a $3$-torsion, we have
{\allowdisplaybreaks
\begin{align*}
& \bigg\{ 1 + 2t^{3e-3}\ep_3y^2z,1+t^{e'}\dfrac{\,\ep_2 y\,}{z} \bigg\}_{\phantom{|}}
 & \\
&\equiv \big\{1+t^{3e'-3}\ep_2\ep_3 y^3,-t^{3e-3}\ep_3y^2z \big\}_{\phantom{|}} \mod U^{4e'-3}H
 & \hbox{(Sublemma \ref{sublem1.thm.cubic}\,(2),\,(3))} \\
&= \big\{ 1 + (\zeta_3-1)^3t^{-3}y^3(1+(\zeta_3-1))^2,\ep_3y^2z \big\}_{\phantom{|}}
 & \hbox{($t^{3e'}\ep_2\ep_3=(\zeta_3-1)^3\zeta_3^2$)} \\
& \equiv \big\{ 1 +(\zeta_3-1)^3t^{-3}y^3,\ep_3y^2z \big\}_{\phantom{|}} \mod U^{4e'-3}H
 & \hbox{(Sublemma \ref{sublem1.thm.cubic}\,(2))} \\
&= -\bigg\{ \dfrac{(1 +(\zeta_3-1)t^{-1}y)^3}{1 +(\zeta_3-1)^3t^{-3}y^3},\ep_3y^2z \bigg\}_{\phantom{|}}
  \\
&= -\bigg\{ 1+\dfrac{3(\zeta_3-1)t^{-1}(y+(\zeta_3-1)t^{-1}y^2)}{1 +(\zeta_3-1)^3t^{-3}y^3},\ep_3y^2z \bigg\}_{\phantom{|}}.
\end{align*}}
The last term belongs to $U^{3e'-1}H$, which proves the assertion.
\par\indent (3) Similarly as for (2), we have
{\allowdisplaybreaks
\begin{align*}
&\bigg\{ 1 + 2t^{3e-3}\ep_3yz^2,1+t^{e'}\dfrac{\,\ep_2 y\,}{z} \bigg\}_{\phantom{|}} & \\
& \equiv \big\{1+t^{3e'-3}\ep_2\ep_3 y^2z,-t^{3e-3}\ep_3yz^2 \big\}_{\phantom{|}} \mod U^{4e'-3}H & \hbox{(Sublemma \ref{sublem1.thm.cubic}\,(2),\,(3))} \\
& = \big\{ 1 + t^{3e'-3}\ep_2\ep_3y^2z, t^{5e'-6}y^3z^3\ep_2\ep_3 \big\}_{\phantom{|}} & \hbox{(Steinberg relation)} \\
& = -\big\{ 1 + t^{3e'-3}\ep_2\ep_3y^2z, (1-\zeta_3)\zeta_3^2 \big\}_{\phantom{|}} &
 \hbox{($t^{5e'}\ep_2\ep_3^2=(1-\zeta_3)^5\zeta_3^4$)} \\
& = -\big\{ 1 + t^{3e'-3}\ep_2\ep_3y^2z, \ep_4(1+\pi^{e-1}u)^{\frac{e}2}t^{e'}\zeta_3^2 \big\}_{\phantom{|}} &  \hbox{($\ep_4:=(1-\zeta_3)\pi^{-\frac e 2}$, \eqref{eq2.0.thm.cubic})} \\
& \equiv -\big\{1 + t^{3e'-3}\ep_2\ep_3y^2z, \ep_4\zeta_3^2 \big\}_{\phantom{|}} \mod U^{5e'-6}H & \hbox{(Sublemma \ref{sublem1.thm.cubic}\,(3))}.
\end{align*}
}The last term is contained in $U^{3e'}H$, because $\ep_4\zeta_3^2$ belongs to $\OO_k^\times$ and $\bF^\times$ is 3-divisible.
This completes the proof of Lemmas \ref{lem.symbol} and \ref{claim3-0.thm.cubic}\,(1).
\end{pf}
\smallskip
We finally prove Lemma \ref{claim3-0.thm.cubic}\,(2).
Put $\ep_3 := 3t^{-3e} \in \Ah^\times$. Since
\addtocounter{equation}{2}
\begin{equation}\label{eq3.3.thm.cubic}
g=\frac{\,x+z\,}{x+y}=\frac{y}{\,z\,}{\cdot}\frac{z(t-y)}{\,y(t-z)\,} = \frac{y}{\,z\,}{\cdot}\bigg(1+\frac{t(z-y)}{\,y(t-z)\,}\bigg),
\end{equation}
we have
{\allowdisplaybreaks \begin{align*}
&\{1+\pi^{e-1}u,g\}
= \bigg\{ 1 + \pi^{e-1}\ep_1(2y^2z+2yz^2+tv),\frac{y}{\,z\,}{\cdot}\bigg(1+\frac{t(z-y)}{\,y(t-z)\,}\bigg)\bigg\} \\
& \quad\equiv \bigg\{ 1 + t^{3e-2}\ep_3v,\frac{y}{\,z\,} \bigg\}
 + \bigg\{ 1 + 2t^{3e-3}\ep_3(y^2z+yz^2),1+\frac{t(z-y)}{\,y(t-z)\,}\bigg\} \\
&\quad \phantom{==} 
 + \bigg\{ 1 + 2t^{3e-3}\ep_3(y^2z+yz^2),\frac{y}{\,z\,}\bigg\}
\end{align*}
}where the congruity holds modulo $U^{3e-1}H$ by Sublemma \ref{sublem1.thm.cubic}\,(1)--(3) and \eqref{eq2.0.thm.cubic}.
One can easily check that the last term belongs to $U^{3e}H$
by similar computations as in Lemma \ref{lem.symbol}\,(3). One can also check that the sum of symbols
\[ \bigg\{ 1 + t^{3e-2}\ep_3v,\frac{y}{\,z\,} \bigg\}
 + \bigg\{ 1 + 2t^{3e-3}\ep_3(y^2z+yz^2),1+\frac{t(z-y)}{\,y(t-z)\,}\bigg\} \in U^{3e-2}H \]
maps to $\ol{\ep_3}(zdy-ydz)$
under the Bloch-Kato isomorphism $\gr_U^{3e-2}H \simeq \Omega^1_\eta$ (with respect to the prime element $t \in \Ah$), by Sublemma \ref{sublem1.thm.cubic}\,(3) and \eqref{eq3.1.thm.cubic}. Here $\ol{\ep_3} \in \kappa(\eta)^\times$ denotes the residue class of $\ep_3=3t^{-3e}$, which belongs to $\bF^\times$ by \eqref{eq2.0.thm.cubic}. This completes the proof of Lemma \ref{claim3-0.thm.cubic}\,(2), Proposition \ref{claim3.thm.cubic} and Theorem \ref{thm.cubic}\,(2).
\end{pf*}

\section{Proof of Theorem \ref{cor0-3}}\label{sect5}
\medskip
In this section, we prove Theorem \ref{cor0-3}. Let $k$ be $p$-adic local field, and let $\OO_k$ be the integer ring of $k$. Put $S:= \Spec(\OO_k)$. Let $\cX$ be a regular scheme which is proper flat of finite type over $S$. Let $Y$ be the closed fiber of $\cX/S$.
\subsection{Proof of Theorem \ref{cor0-3}\,(1)}\label{sect5-1}
Let \[  \langle \;,\; \rangle : \CH_0(\gfiber) \times \Br(\gfiber) \lra \qz \] be the Brauer-Manin pairing \eqref{eq01} in the introduction, and let $\omega$ be an element of $\Br(\gfiber)$ with  $\langle c,\omega \rangle=0$ for any $c \in \CH_0(\gfiber)$.
It suffices to show that $\omega$ belongs to $\Br(\cX)$. Indeed, $\omega$ is $0$-unramified in the sense of Definition \ref{def:ur}\,(2), and the assertion follows from Corollary \ref{cor:key}, where we have assumed the purity of Brauer groups for $\cX$.
\subsection{Proof of Theorem \ref{cor0-3}\,(2)}\label{sect5-2}
We start the proof of Theorem \ref{cor0-3}\,(2), which will be finished in \S\ref{sect5-3} below.
Let us recall the map \eqref{eq03}, which we denote by $\Phi$ in what follows:
\[ \Phi : \A_0(\gfiber) \lra \Hom(\Br(\gfiber)/\Br(k)+\Br(\cX), \Q/\Z) \]
By Theorem \ref{cor0-3}\,(1), $\Phi$ has dense image.
To prove that $\Phi$ is surjective, it is enough to show
\[ \Image(\Phi) \simeq \Z_p^{\oplus r} \oplus T \] for some non-negative integer $r$ and some finite group $T$.
We are thus reduced to the following proposition, where we do not assume the purity of Brauer groups:
\begin{prop}\label{prop:finite}
Put $B:=\Br(\gfiber)/\Br(k)+\Br(\cX)$ and $D_{\ell}:= \Hom(B_{\ltor},\Q/\Z)$ for a prime number $\ell$. Then{\rm:}
\begin{enumerate}
\item[(1)] $D_{\ell}$ is finitely generated over $\Z_{\ell}$ for any $\ell$.
\item[(2)] If $\ell\not=p$, then $D_{\ell}$ is finite and the map $\Phi_{\ell} : \A_0(\gfiber) \ra D_{\ell}$ induced by $\Phi$ is surjective. The image of $\Phi_{p} : \A_0(\gfiber) \ra D_p$ is a $\zp$-submodule of $D_p$.
\item[(3)] $D_{\ell}$ is zero for almost all $\ell \not = p$.
\end{enumerate}
\end{prop}
\begin{rem}
The `$\ell \not = p$' case of Proposition \ref{prop:finite}\,(2) is proved in \cite{cts} Corollaire 2.6.
\end{rem}
\noindent
{\it Proof of Proposition \ref{prop:finite}.}
The assertion (1) is obvious. We prove (2). Note that the map $\Phi_{\ell}$ with $\ell \not = p$ has dense image with respect to the $\ell$-adic topology on $D_{\ell}$ (by the $\ell$-primary part of Theorem \ref{cor0-3}\,(1) and the absolute purity).
By a standard norm argument, we may suppose that $\gfiber$ has a $k$-rational point. Then we have a surjective map
\[ \bigoplus_{C} \ \A_0(C) \lra \A_0(\gfiber). \]
Here $C$ ranges over the smooth integral curves over $k$ which are finite over $\gfiber$.
Since $k$ is a $p$-adic field by assumption, we have \[ \A_0(C) \simeq \Z_p^{\oplus r_C} \oplus T_C \] for a non-negative integer $r_C$ and a finite group $T_C$ by a theorem of Mattuck \cite{mat}. In case $\ell \not= p$, these facts and (1) imply that $\Image(\Phi_{\ell})$ is finite and dense, so that $\Phi_{\ell}$ is surjective.
To prove the assertion for $\Phi_p$, we need to show that the composite map
\[\xymatrix{ \Z_p^{\oplus r_C} \; \ar@<-1pt>@{^{(}->}[r] & \A_0(C) \ar[r] & \A_0(\gfiber) \ar[r] & D_p }\]
is a homomorphism of $\Z_p$-modules, that is, continuous with respect to the $p$-adic topology. Let $f_C: \Br(k) \ra \Br(C)$ be the natural restriction map. Since the above composite map factors through the map
\[\xymatrix{ \Z_p^{\oplus r_C} \; \ar@<-1pt>@{^{(}->}[r] & \A_0(C) \ar[r] & \Hom(\Coker(f_C)_{\ptor},\qz), }\]
it suffices to see the continuity of this map, which is a consequence of \cite{ss:localring} Theorem (9.2).
Thus we obtain Proposition \ref{prop:finite}\,(2). \par
We next show Proposition \ref{prop:finite}\,(3). Suppose $\ell \not =p$. Since we have $\Image(\Phi_{\ell})=D_{\ell}$ and
\[ \Image(\Phi_{\ell}) \simeq \Image(\A_0(\gfiber) \ra \H^{2N}(\gfiber,\zl(N))) \]
with $N:=\dim(\gfiber)$, the problem is reduced to the case where $\cX$ has strict semistable reduction over $S$ by the alteration theorem of de Jong \cite{dJ} and a standard norm argument using the functoriality of cycle class maps, where `strict' means that all irreducible components of $Y$ are smooth.
We prove that $\bigoplus_{\ell \not = p}\ D_{\ell}$ is finite, assuming that $\cX/S$ has strict semistable reduction.
For a torsion abelian group $M$, let $M'$ be its prime-to-$p$ part. For $n \in \Z$, let $\qzb(n)$ be the \'etale sheaf $\bigoplus_{\ell \not =p}\ \qzl(n)$. Consider a commutative diagram with exact rows
\[ \xymatrix{ 0 \ar[r] & \Pic(\cX) \otimes \qzb \ar[r] \ar@{->>}[d] & \H^2(\cX,\qzb(1)) \ar[r] \ar[d]_a & \Br(\cX)' \ar[r] \ar@<-1pt>@{^{(}->}[d] & 0 \\ 0 \ar[r] & \Pic(\gfiber) \otimes \qzb \ar[r] & \H^2(\gfiber,\qzb(1)) \ar[r] & \Br(\gfiber)' \ar[r] & 0. }\]
Since the left vertical arrow is surjective, we have
\[ \Br(\gfiber)'/\Br(\cX)' \simeq \Coker(a) \simeq \ker\big( \H^3_Y(\cX,\qzb(1)) \ra \H^3(\cX,\qzb(1)) \big). \]
Our task is to show that the complex
\addtocounter{equation}{2}
\begin{equation}\label{BM.1.1}
  \Br(k)' \lra \H^3_Y(\cX,\qzb(1)) \lra \H^3(\cX,\qzb(1))
\end{equation} has finite cohomology group at the middle.
Let $\knr$ be the maximal unramified extension of $k$, and put $\Xnr:=\cX \otimes_{\OO_k}\OO_{\knr}$ and $\Ynr :=Y \otimes_{\rf} \ol {\rf}$. In view of the short exact sequences
\begin{align*}
& 0 \to \H^1(\rf, \H^2_{\ol Y}(\Xnr,\qzb(1))) \to \H^3_Y(\cX,\qzb(1)) \to \H^3_{\ol Y}(\Xnr,\qzb(1))^{G_{\rf}} \to 0,\\
& 0 \to \H^1(\rf, \H^2(\Xnr,\qzb(1))) \to \H^3(\cX,\qzb(1)) \to \H^3(\Xnr,\qzb(1))^{G_{\rf}} \to 0
\end{align*}
obtained from Hochschild-Serre spectral sequences and the fact $\cd(\rf)=1$, we are reduced to the following lemma:
\stepcounter{thm}
\begin{lem}\label{lem:finite}
\begin{enumerate}
\item[(1)] The group $\H^3_{\ol Y}(\Xnr,\qzb(1))^{G_{\rf}}$ is finite.
\item[(2)] The composite map $\Br(k)' \to \H^3_Y(\cX,\qzb(1)) \to \H^3_{\ol Y}(\Xnr,\qzb(1))^{G_{\rf}}$ is zero.
Consequently, the complex \eqref{BM.1.1} induces a complex
\addtocounter{equation}{1}
\begin{equation}\label{BM.1.2}
  \Br(k)' \lra \H^1(\rf, \H^2_{\ol Y}(\Xnr,\qzb(1))) \lra \H^1(\rf, \H^2(\Xnr,\qzb(1))),
\end{equation}
\item[(3)] The complex \eqref{BM.1.2} has finite cohomology group at the middle.
\end{enumerate}
\end{lem}
\noindent
We first show Lemma \ref{lem:finite}\,(1).
For $q \geq 1$, let $\ol Y {}^{(q)}$ be the disjoint union of the intersections of $q$ distinct irreducible components of $\ol Y$.
Note that all connected components of $\ol Y {}^{(q)}$ are smooth proper varieties over $\ol \rf$ of dimension $\dim(\cX)-q$.
By the Mayer-Vietoris spectral sequence
\begin{equation}\label{BM.1.MV}
 E_1^{u,v}=\H^{2u+v-2}(\ol Y{}^{(-u+1)},\qzb(u)) \Lra \H^{u+v}_{\ol Y}(\Xnr,\qzb(1))
\end{equation}
(cf.\ \cite{RZ} Satz 2.21), we have the following exact sequence of $G_{\rf}$-modules:
\[ 0 \lra \H^1(\ol Y{}^{(1)},\qzb) \lra \H^{3}_{\ol Y}(\Xnr,\qzb(1)) \lra \H^0(\ol Y{}^{(2)},\qzb(-1)). \]
Hence the assertion follows from the finiteness of $\H^1(\ol Y{}^{(1)},\qzb)^{G_{\rf}}$ due to Katz-Lang \cite{KL}.
Lemma \ref{lem:finite}\,(2) is a consequence of (1) and the fact that $\Br(k)'$ is divisible.
\par
We next show Lemma \ref{lem:finite}\,(3). For an abelian group $M$, let $M_{\L\text{-}\div}$ be the subgroup of elements which are divisible in $M$ by all integers prime to $p$, that is,
\[ M_{\L\text{-}\div} := \bigcap_{n \in \Z,\; (n,p)=1}~n \cdot M. \]
We define a discrete $G_{\rf}$-module $\Xi$ as $\Pic(\Xnr)/\Pic(\Xnr)_{\L\text{-}\div}$. For this group we will prove:
\addtocounter{thm}{2}
\begin{lem}\label{lem:finite2}
\begin{enumerate}
\item[(1)] $\Xi$ is a finitely generated abelian group.
\item[(2)] For any positive integer $n$ prime to $p$, ${}_n\Pic(\Xnr)$ is finite. Consequently, $\Pic(\Xnr)_{\L\text{-}\div}$ is divisible by integers prime to $p$ {\rm(}cf.\ \cite{Ja1} {\rm\S4)}, and we have a natural injective map
\[\xymatrix{ \alpha : \Xi \otimes \qzb \; \ar@<-1pt>@{^{(}->}[r] & \H^2(\Xnr, \qzb(1)). }\]
\item[(3)]
The map \[ \alpha' : \H^1(\rf, \Xi \otimes \qzb) \lra \H^1(\rf,\H^2(\Xnr, \qzb(1))) \] induced by $\alpha$ has finite kernel.
\end{enumerate}
\end{lem}
\noindent
The proof of this lemma will be given in \S\ref{sect5-3} below. We finish our proof of Lemma \ref{lem:finite}\,(3) (and Proposition \ref{prop:finite}\,(3)), admitting this lemma. By Lemma \ref{lem:finite2}\,(1), we have the following sequence of finitely generated abelian groups with equivariant $G_{\rf}$-action:
\addtocounter{equation}{1}
\begin{equation}\label{BM.1.3}
\xymatrix{ \us{\phantom{\ol {(^(}}}{\Z} \ar[rr]^-{\text{diagonal}} && \displaystyle \bigoplus_{y \in (\ol Y)^0} \ \Z \ar[r]^-g & \us{\phantom{\ol {(^c}}}{\Xi ,} } \end{equation}
where $g$ is induced by the Gysin map $\bigoplus_{y \in (\ol Y)^0}\ \Z \ra \Pic(\Xnr)$.
This sequence is a complex by the semistability of $\cX/S$.
Because the $\ell$-primary part of the complex \eqref{BM.1.2} has finite cohomology group for any $\ell \not =p$
 by Proposition \ref{prop:finite}\,(2), one can easily check that the complex \eqref{BM.1.3} has finite cohomology group as well
 (in fact, \eqref{BM.1.3} is exact because the cokernel of the first diagonal map is torsion-free).
Hence the cohomology group of the induced complex
\begin{equation}\label{BM.1.4}
\xymatrix{ \us{\phantom{(^i}}{\qz} \ar[rr]^-{\text{diagonal}} && \displaystyle \bigoplus_{y \in (\ol Y)^0} \ \qz \ar[rr]^-{g\otimes \qz} && \us{\phantom{\ol (}}{\Xi \otimes \qz,} } 
\end{equation} is finite and its order is the same as that of
\[ \Coker\big( g: \textstyle \bigoplus_{y \in (\ol Y)^0} \ \Z \lra  \Xi/ \Xi_{\tor} \big){}_{\tor}. \]
We denote its order by $b$. Now we show that the complex \eqref{BM.1.2} has finite cohomology group.
By a standard norm argument, we may suppose that $G_{\rf}$ acts trivially on the groups in \eqref{BM.1.3}.
Then consider a commutative diagram of complexes
\begin{equation}\label{BM.1.5}
\hspace{-10pt}\xymatrix{ \H^1(\rf, \qzb) \ar[r] & \bigoplus_{y \in Y^0} \ \H^1(\rf,\qzb) \ar[r] \ar[d]_{\wr\hspace{-2pt}} & \H^1(\rf,\Xi \otimes \qzb) \ar[d]_{\alpha'} \\
\Br(k)' \ar[r] \ar[u]^{\wr\hspace{-2pt}} & \H^1(\rf, \H^2_{\ol Y}(\Xnr,\qzb(1))) \ar[r] & \H^1(\rf, \H^2(\Xnr,\qzb(1))),}
\end{equation}
where the upper row is induced by \eqref{BM.1.4}, the lower row is the complex \eqref{BM.1.2} in question
 and the bijectivity of the central vertical arrow is obtained from the spectral sequence \eqref{BM.1.MV};
 the commutativity of the left square follows from the semistability of $\cX/S$.
The upper row has finite cohomology group of order $b$ up to a power of $p$, and the right vertical arrow $\alpha'$ has finite kernel by Lemma \ref{lem:finite2}\,(3). Hence the lower row has finite cohomology group of order dividing $b \cdot \#(\ker(\alpha'))$. This completes the proof of Lemma \ref{lem:finite}, Proposition \ref{prop:finite} and Theorem \ref{cor0-3}\,(2), assuming Lemma \ref{lem:finite2}. \hfill $\square$
\subsection{Proof of Lemma \ref{lem:finite2}}\label{sect5-3}
We first prove (1). We change the notation slightly, and put $\Xi(Z):=\Pic(Z)/\Pic(Z)_{\L\text{-}\div}$ for a scheme $Z$.
For a smooth variety $Z$ over a field, let $\NS(Z)$ be the N\'eron-Severi group of $Z$.
The natural map $\Xi(\Xnr) \ra \Xi(\Ynr)$ is injective by the proper base-change theorem:
\[ \xymatrix{ \Xi(\Xnr) \ar[r] \ar@{^{(}->}[d] & \Xi(\Ynr) \ar@{^{(}->}@<.5mm>[d] \\
 \displaystyle  \prod_{\ell \not =p} \ \H^2(\Xnr,\zl(1)) \ar[r] \ar@{}@<-.7mm>[r]^\sim & \displaystyle  \prod_{\ell \not =p} \ \H^2(\ol Y,\zl(1)).} \]
We show that $\Xi(\Ynr)$ is finitely generated. Since $Y$ has simple normal crossings on $\cX$, we have the following exact sequence of sheaves on $Y_{\et}$:
\[ 0 \lra \cO_Y^\times \lra \Gm,_{Y^{(1)}} \lra \Gm,_{Y^{(2)}} \lra \dotsb, \]
where $Y^{(1)}$ ($q \geq 1$) is as in the proof of Lemma \ref{lem:finite}\,(1).
By this exact sequence, it is easy to see that the kernel of the natural map $\Pic(\ol Y) \ra \Pic(\ol Y{}^{(1)})$ is an extension of a finite group by a torsion divisible group. Hence the induced map
\[ \Xi(\ol Y) \lra \Xi(\ol Y{}^{(1)}) \simeq \NS(\ol Y{}^{(1)})/\NS(\ol Y{}^{(1)})_{\ptor} \]
 has finite kernel, and the assertion follows from the fact that the last group is a finitely generated abelian group.
Thus we obtain Lemma \ref{lem:finite2}\,(1).
\par
Lemma \ref{lem:finite2}\,(2) follows from the finiteness of the groups $\H^1(\Xnr,\mu_n)$ with $(n,p)=1$.
The details are straight-forward and left to the reader.
\par
We prove Lemma \ref{lem:finite2}\,(3). By the proof of (1), the natural map $\Xi(\Xnr) \ra \Xi(\ol Y{}^{(1)})$ has finite kernel and
 $\Xi(\ol Y{}^{(1)})$ is finitely generated, which implies that the left vertical arrow in the following commutative diagram
 has finite kernel by a standard norm argument:
\[\xymatrix{ \H^1(\rf, \Xi(\Xnr) \otimes \qzb) \ar[r]^-{\alpha'} \ar[d] & \H^1(\rf,\H^2(\Xnr, \qzb(1))) \ar[d] \\
 \H^1(\rf, \Xi(\ol Y{}^{(1)}) \otimes \qzb) \ar[r] & \H^1(\rf,\H^2(\ol Y{}^{(1)}, \qzb(1))), }\]
where the bottom horizontal arrow is defined in the same way as $\alpha'$ and it has finite kernel by \cite{SaSa} Lemma 6.7
 (note that $\Xi(\ol Y{}^{(1)}) = \NS(\ol Y{}^{(1)})/\NS(\ol Y{}^{(1)})_{\ptor}$). Hence $\alpha'$ has finite kernel as well.
This completes the proof of Lemma \ref{lem:finite2}.
\qed
\subsection{Appendix: Degree of 0-cycles}\label{sect5-4}
Let $\cX$ be a regular scheme of finite type over $S=\Spec(\OO_k)$. We do not assume the properness of $\cX/S$ here.
Put $Q := \gfiber_0$, and let $Q_{\fin}$ be the subset of $Q$ consisting of all closed points on $\gfiber$ whose closure in $\cX$ are finite over $S$. If $\cX$ is proper over $S$ then we have $Q_{\fin}=Q$.
For $y \in Y^0$, let $e_y$ be the multiplicity of $y$ in $\cX \otimes_{\OO_k}\bF$, and let $f_y$ be the degree over $\bF$
 of the algebraic closure of $\bF$ in $\kappa(y)$.
\begin{thm}\label{thm:degree}
Assume that the purity of Brauer groups holds for $\cX$. Then the following three numbers are equal to one another$:$
\begin{enumerate}
\item[$N_1$] {\rm:} the order of the kernel of the composite map $\Br(k) \ra \Br(\gfiber) \ra \Br(\gfiber)/\Br(\cX).$
\item[$N_2$] {\rm:} the greatest common measure of the degrees $[\kappa(v):k]$ with $v \in Q_{\fin}$.
\item[$N_3$] {\rm:} the greatest common measure of the integers $e_y \cdot f_y$ with $y \in Y^0$.
\end{enumerate}
\end{thm}
\begin{rem}
This result was proved by Colliot-Th\'el\`ene and Saito in \cite{cts} Th\'eor\`eme 3.1, up to powers of $p$.
However, their arguments work including powers of $p$ after some modification, which we show in what follows.
\end{rem}
\begin{pf}
By the arguments in \cite{cts} Th\'eor\`eme 3.1, we have $N_1 \vert N_2 \vert N_3$.
Our task is to show $N_1=N_3$. Let us consider the composite map
\[ \alpha : \Br(k) \lra \Br(\gfiber)/\Br(\cX) \os{\beta}\lra \bigoplus_{y \in Y^0}\ \H^3_{y}(\cX,\Gm), \]
where the second map $\beta$ denotes the composite map
\[ \xymatrix{ \us{\phantom{(}}{\beta : \Br(\gfiber)/\Br(\cX)} \; \ar@<-1pt>@{^{(}->}[r] & \us{\phantom{(}}{\H^3_{Y}(\cX,\Gm)} \ar[r] & \displaystyle  \bigoplus_{y \in Y^0}\ \H^3_{y}(\cX,\Gm), }\] and the last map is injective by the purity assumption on Brauer groups (cf.\ Remark \ref{rem:ur}\,(3)).
Hence it suffices to show the following lemma:
\begin{lem}\label{lem:degree}
Let $\eta$ be a generic point of $Y$. Then the kernel of the composite map
\[\xymatrix{ \us{\phantom{(}}{\alpha_{\eta} : \Br(k)} \;\ar@<-1pt>@{^{(}->}[r]^-{\alpha} & \displaystyle \bigoplus_{y \in Y^0}\ \H^3_{y}(\cX,\Gm) \ar[r] & \us{\phantom{1}}{\H^3_{\eta}(\cX,\Gm)} }\] is isomorphic to $\Z/e_{\eta}f_{\eta}\Z$.
\end{lem}
\medbreak
\noindent
{\it Proof of Lemma \ref{lem:degree}.}
Let $\Ah$ (resp.\ $A_{\ol \eta}$) be the henselization (resp.\ strict henselization) of $\cO_{\sscX,\eta}$ at its maximal ideal,
  and let $G_{\eta}$ be the absolute Galois group of $\kappa(\eta)$.
By the Hochschild-Serre spectral sequence
\[ E_2^{u,v}=\H^u_{\Gal}(G_{\eta},\H^v_{\ol \eta}(A_{\ol \eta},\Gm)) \Lra \H^{u+v}_{\eta}(A_{\eta},\Gm) (= \H^{u+v}_{\eta}(\cX,\Gm)), \]
and a standard purity for $\Gm$ (\cite{gr} III.6), we have the following exact sequence:
\[ 0 \lra \H^1(\eta,\Q/\Z) \lra \H^3_{\eta}(\Ah,\Gm) \lra \H^3_{\ol \eta}(A_{\ol \eta},\Gm)^{G_{\eta}}. \]
Next let $s$ be the closed point of $S$, and let $\knr$ be the maximal unramified extension of $k$.
The same computation for $S$ yields an isomorphism $\H^1(\bF,\Q/\Z) \simeq \H^3_{s}(S,\Gm) (\simeq  \Br(k))$,
 because we have $\H^3_{\ol s}(\O_{S,\ol s}^{\sh},\Gm)\simeq \Br(\knr)=0$ (cf.\ \cite{Se} II.3.3).
Thus we obtain the following commutative diagram with exact rows:
\[ \xymatrix{ 0 \ar[r] & \H^1(\bF,\Q/\Z) \ar[r] \ar@<-.5mm>@{}[r]^\sim \ar[d]_{\alpha_{\eta}'} & \H^3_{s}(S,\Gm) \ar[r] \ar[d]_{\alpha_{\eta}} & 0 \ar[d] \\
0  \ar[r] & \H^1(\eta,\Q/\Z) \ar[r] & \H^3_{\eta}(\Ah,\Gm)\ar[r] & \H^3_{\ol \eta}(A_{\ol \eta},\Gm)^{G_{\eta}},}\]
where $\alpha_{\eta}'$ denotes the map induced by the right square.
One can easily show that \[ \alpha_{\eta}'=e_{\eta} \cdot \Res, \] where $\Res$ denotes the natural restriction map
 $\H^1(\bF,\Q/\Z) \to \H^1(\eta,\Q/\Z)$ induced by the structure map $\eta \ra \Spec(\bF)$.
Thus we have \[ \ker(\alpha_{\eta}) \simeq \ker(\alpha_{\eta}') \simeq \Z/e_{\eta}f_{\eta}\Z. \]
This completes the proof of Lemma \ref{lem:degree} and Theorem \ref{thm:degree}.
\end{pf}

\section{Proof of Theorem \ref{thm:surj}}\label{sect6}
\medskip
In this section, we prove Theorem \ref{thm:surj}. Let the notation be as in \S\ref{sect1-2}.
Put $S:= \Spec(\OO_k)$ and let $Y$ be the divisor on $\cX$ defined by the radical of $(p) \subset \cO_{\sscX}$.
\subsection{Key diagram}\label{sect6-1}
For $v \in \gfiber_0$, the closure $\ol {\{ v \} } \subset \cX$ contains exactly one closed point of $Y$ by the properness of $\cX$.
\begin{defn}
We define a map of sets $\sp: \gfiber_0 \ra Y_0$ by the law that $\sp(v)=x$ if and only if $x$ is the closed point of
 $\ol{\{ v \} } \subset \cX$. Since $\OO_k$ is henselian local, there is a natural identification of sets
\[ \sp^{-1}(\{ x \}) = \{ \hbox{closed points on } \Spec(\cO_{\sscX,x}^{\h}\ip)\} \qquad (x \in Y_0) \]
where $\cO_{\sscX,x}^{\h}$ denotes the henselization of $\cO_{\sscX,x}$ at its maximal ideal.
\end{defn}
\noindent
We construct the key diagram \eqref{CD1} below. Let $U \subset Y$ be a regular dense open subset.
For $x \in U_0$, put $A_x:=\cO_{\sscX,x}^{\h}$ and let $\psi_{x,p^r}$ be the composite map
\[ \xymatrix{ \us{\phantom{(}}{\psi_{x,p^r} : {}_{p^r}\Br(A_x \ip )} \ar[rr]^{\text{specialization}} &&
   \displaystyle \prod_{v \in \sp^{-1}(\{ x \})} \ {}_{p^r}\Br(v)
 \ar[rr]^{\text{invariant}}_\sim && \displaystyle \prod_{v \in \sp^{-1}(\{ x \})} \ \Z/p^r\Z. } \]
We put $Q_U:=\sp^{-1}( U_0 )$, and define a map
\[ \psi_{p^r} : \prod_{x \in U_0} \ {}_{p^r}\Br(A_x \ip ) \lra \prod_{v \in Q_U} \ \Z/p^r\Z \]
as the direct product of $\psi_{x,p^r}$'s with $x \in U_0$. Next we construct a canonical map
\[ \theta_{p^r} : \H^3_{Y}(\cX,\T_r(1)) \lra \prod_{x \in U_0} \ {}_{p^r}\Br(A_x \ip), \]
where $\T_r(1)=\T_r(1)_{\scX}$ denotes the $p$-adic \'etale Tate twist \cite{Sat2}. For this we will prove
\begin{lem}\label{lem9-1}
Let $x$ be a closed point on $U$ and let $Y_x$ be the divisor on $\cX_x:=\Spec(A_x)$ defined by the radical of $(p) \subset A_x$.
Then there is a canonical isomorphism \[ {}_{p^r}\Br(A_x\ip) \simeq \H^3_{Y_x}(\cX_x,\T_r(1)). \]
\end{lem}
\noindent
By this lemma, we define the desired map $\theta_{p^r}$ as the restriction map
\[ \H^3_{Y}(\cX,\T_r(1)) \lra \prod_{x \in U_0} \ \H^3_{Y_x}(\cX_x,\T_r(1)) \simeq  \prod_{x \in U_0} \ {}_{p^r}\Br(A_x \ip ). \]
\begin{pf*}{\it Proof of Lemma \ref{lem9-1}}
Consider the localization exact sequence
\[ \H^2(\cX_x,\T_r(1)) \lra \H^2(\cX_x\ip,\mu_{p^r}) \lra \H^3_{Y_x}(\cX_x,\T_r(1)) \lra \H^3(\cX_x,\T_r(1)), \]
where we have used the property that $\T_r(1)$ is isomorphic to $\mu_{p^r}$ outside of characteristic $p$ (cf.\ \cite{Sat2}).
Since $A_x \ip$ is a unique factorization domain, we have $\H^2(\cX_x\ip,\mu_{p^r}) \simeq {}_{p^r}\Br(A_x \ip)$.
On the other hand, we have
\begin{align*}
\H^2(\cX_x,\T_r(1)) &\os{(*)}{\simeq} {}_{p^r}\Br(\cX_x) \simeq {}_{p^r}\Br(x) = 0, \\
\H^3(\cX_x,\T_r(1)) & \simeq  \H^3(x, \T_r(1)\vert_x) = 0,
\end{align*}
where the isomorphism $(*)$ follows from the Kummer theory for $\cO_{\sscX}^\times$ (cf.\ \cite{Sat2} \S4.5),
 and the last vanishing follows from the facts that $\T_r(1)$ is concentrated in degree $0$ and $1$ (cf.\ loc.\ cit.) and that $\cd(x)=1$.
The assertion follows from these facts.
\end{pf*}
\begin{cor}\label{sublem1}
Under the notation in Lemma \ref{lem9-1}, the Gysin map {\rm(}\cite{Sat2} Theorem {\rm 4.4.7)}
\[ \H^1(Y_x,\Z/p^r\Z ) \lra \H^3_{Y_x}(\cX_x,\T_r(1)). \]
is injective.
\end{cor}
\begin{pf}
By Hensel's lemma, one can take a one-dimensional regular closed subscheme $Z \subset \cX_x$ which is \'etale over $S$. Let $y$ be the generic point of $Z$. Then the composite map
\[ \Z/p^r\Z \simeq \H^1(Y_x,\Z/p^r\Z) \to \H^3_{Y_x}(\cX_x,\T_r(1)) \simeq {}_{p^r}\Br(A_x\ip) \to {}_{p^r}\Br(y) \simeq \Z/p^r\Z \]
is bijective by the construction of the isomorphism of Lemma \ref{lem9-1}.
\end{pf}

We return to the proof of Theorem \ref{thm:surj}, and define $\alpha$ as the composite map
\[\xymatrix{ \alpha: \CH^{d-1}(\cX)/p^r \ar[rr]^-{\varrho_{p^r}^{d-1}} && \H^{2d-2}(\cX,\T_r(d-1))
   \isom \H^3_Y(\cX,\T_r(1))^*, }\]
where for a $\Z/p^r\Z$-module $M$, $M^*$ denotes $\Hom(M,\Z/p^r\Z)$,
  and the last isomorphism follows from the arithmetic duality (\cite{Sat2} Theorem 10.1.1).
To prove Theorem \ref{thm:surj}, it remains to show the following:
\begin{lem}\label{lem:CD}
\begin{enumerate}
\item[(1)]
The following diagram is commutative{\rm:}
\addtocounter{equation}{3}
\begin{equation}\label{CD1}
\xymatrix{ \CH^{d-1}(\cX)/p^r \ar[rr]^-{\alpha} && \H^3_Y(\cX,\T_r(1))^* \\
 \displaystyle  \bigoplus_{v \in Q_U}\  \Z/p^r\Z \ar[rr]^-{(\psi_{p^r})^*} \ar[u] &&
 \displaystyle \bigoplus_{x \in U_0}\ {}_{p^r}\Br(A_x \ip)^*, \ar[u]_{(\theta_{p^r})^*}}
\end{equation}
where the left vertical arrow sends $v \in Q_U$ to the class of $v$.
\item[(2)]
The maps $\theta_{p^r}$ and $\psi_{p^r}$ are injective. Consequently, $(\theta_{p^r})^*$ and $(\psi_{p^r})^*$ are surjective.
\end{enumerate}
\end{lem}
\noindent
We prove (2) in \S\ref{sect6-2}, and then prove (1) in \S\ref{sect6-3} below.
\subsection{Proof of Lemma \ref{lem:CD}\,(2)}\label{sect6-2}
Put $Z:=Y \sm U$ and $\cX' := \cX \sm Z$. For $x \in U_0$, let $\cX_x$ and $Y_x$ be as in Lemma \ref{lem9-1}. The map $\psi_{p^r}$ is injective by Corollary \ref{thm:key}. Indeed this injectivity is immediately reduced to the case $\zeta_p \in k$ by a standard norm argument, because $\cX'$ is smooth over $S$ by the semistability assumption on $\cX$. To show the injectivity of $\theta_{p^r}$, we consider the following commutative diagram:
\[ \xymatrix{ \theta_{p^r} : \H^3_{Y}(\cX,\T_r(1)) \; \ar@<-1pt>@{^{(}->}[r]^-{a} & \H^3_{U}(\cX',\T_r(1)) \ar[r] & \prod_{x \in U_0} \ \H^3_{Y_x}(\cX_x,\T_r(1)) \\ & \H^1(U,\Z/p^r\Z) \; \ar@<-1pt>@{^{(}->}[r]^-{b} \ar[u]_c & \prod_{x \in U_0} \ \H^1(Y_x,\Z/p^r\Z) \ar@{^{(}->}[u]_{c'}, } \]
where $a$ and $b$ are restriction maps, and $c$ and $c'$ are Gysin maps (\cite{Sat2} Theorem 4.4.7).
The arrow $a$ is injective by the purity for $\T_r(1)$ (cf.\ loc.\ cit.), and $b$ is injective by the \v{C}ebotarev density theorem \cite{Se2} Theorem 7.
The arrow $c'$ is injective by Lemma \ref{lem9-1} and Corollary \ref{sublem1}.
It remains to show that the induced map $\Coker (c) \ra \Coker (c')$ is injective.
Let $i$ be the closed immersion $U \hra \cX'$. Since we have $\Coker (c) \subset \H^0(U,R^3i^!\T_r(1))$ and
\[ \Coker (c') \subset \prod_{x \in U_0} \ \H^0(Y_x,(R^3i^!\T_r(1))\vert_{Y_x}), \]
it suffices to show that for an \'etale sheaf $\FF$ on a noetherian scheme $W$, the restriction map
\[ \H^0(W,\FF)  \lra \prod_{x \in W_0} \ \H^0\big(\Spec(\O_{W,x}^{\h}),\FF \vert_{\Spec(\O_{W,x}^{\h})}\big) \]
is injective, which follows from a standard argument using the induction on $\dim(W)$.
Thus we obtain Lemma \ref{lem:CD}\,(2).
\subsection{Proof of Lemma \ref{lem:CD}\,(1)}\label{sect6-3}
Put $\Lam:= \Z/p^r\Z$ and $N:=d-1$ for simplicity. By the definitions of the maps in \eqref{CD1}, it is enough to show that the following diagram commutes for each $v \in Q_U$:
\[\xymatrix{
\quad\;\;\; \H^0(C_v,\Lam) \ar@<-42pt>[d]_{f_*} \quad\; \times \H^3_{x}(C_v,\T_r(1)_{C_v}) \ar[r] \ar@{}[rd]|{\qquad\qquad(1)} & \H^3_{x}(C_v,\T_r(1)_{C_v}) \ar[r]^-{\tr_{(C_v,x)}} \ar[d]_{f_*} \ar@{}[rd]|{\quad(2)} & \Lam \ar@{=}[d] \\
\H^{2N}(\cX,\T_r(N)_{\scX}) \times \H^3_{Y}(\cX,\T_r(1)_{\scX}) \ar[r] \ar@<-54pt>[u]_{f^*} & \H^{2d+1}_{Y}(\cX,\T_r(d)_{\scX}) \ar[r]^-{\tr_{(\cX,Y)}} & \Lam, \hspace{-2pt} } \]
where $C_v$ denotes the normalization of $\ol { \{ v \} } \subset \cX$, $x$ denotes the closed point of $C_v$ and $f$ denotes the canonical finite morphism $C_v \ra \cX$. See \cite{Sat2} Theorem 10.1.1 for the trace maps.
The arrows $f_*$ arise from the following relative trace morphism with $n=0,1$ (loc.\ cit.\ Theorem 7.1.1):
\[ \gys_f: Rf_*\T_r(n)_{C_v} \lra \T_r(n+N)_{\scX}\,[2N] \quad \hbox{ in } \; D^b(\cX_{\et},\Lam), \]
and the arrow $f^*$ arises from the pull-back morphism (loc.\ cit.\ Proposition 4.2.8)
\[ \res^f: \T_r(1)_{\scX} \lra Rf_*\T_r(1)_{C_v} \quad \hbox{ in } \; D^b(\cX_{\et},\Lam). \]
The commutativity of the square (2) follows from a similar argument as for \cite{Sat2} Lemma 10.2.1.
We prove the commutativity (1) of pairings. Consider the following commutative diagram:
\[\xymatrix{ Y \;\ar@<-1pt>@{^{(}->}[r]^i \ar[d]_\gamma \ar@{}[rd]|{\square} & \cX \ar[d]_g & \ar[l]_f \ar[ld]^{\pi} C_v \\ s  \;\ar@<-1pt>@{^{(}->}[r]^h & S, }\]
where $s$ denotes the closed point of $S$. Note that we have the base-change isomorphism $R\gamma_*Ri^! = Rh^!Rg_*$ by Deligne \cite{sga4} XVIII.3.1.12.3.
To prove the commutativity of pairings in question, it suffices to show that the following diagram commutes in $D^-(s_{\et},\Lam)$:
\begin{equation}\label{CD4}
\hspace{-65pt} \xymatrix{ h^* R\pi_* \Lam_{C_v} \otimes^{\L} Rh^!Rg_*\T_r(1)_{\scX} \ar[rr]^-{\id \otimes^{\L} \, R h^!R g_* (\res^f)}  \ar[dd]_{h^* Rg_*(\gys_f) \otimes^{\L} \, \id } && h^* R\pi_* \Lam_{C_v} \otimes^{\L} Rh^! R\pi_* \T_r(1)_{C_v}
 \ar@<-8pt>[d]^{{\text{product}}} \\ &&  Rh^!R\pi_*\T_r(1)_{C_v} \ar@<-8pt>[d]^{R h^!R g_*(\gys_f)} \\
h^* Rg_* \T_r(N)_{\scX}\,[2N] \otimes^{\L} Rh^!Rg_*\T_r(1)_{\scX} \ar[rr]^-{\text{product}} && Rh^!Rg_*\T_r(d)_{\scX}\,[2N].  } \hspace{-80pt} \end{equation}
Here the arrows `product' are induced by the canonical morphism
\[ h^* \cK \otimes^{\L} Rh^! \cL \lra Rh^!(\cK \otimes^{\L} \cL) \qquad \hbox{($\cK,\cL \in D^b(S_{\et},\Lam)$)} \]
and the product structure of Tate twists (\cite{Sat2} Proposition 4.2.6). We defined $Rh^!$ for unbounded objects using a result of Spaltenstein \cite{Spa} Theorem A. Finally one can easily check the commutativity of \eqref{CD4} by applying $Rh^!Rg_*$ to the following commutative diagram in \cite{Sat2} Corollary 7.2.4:
\begin{equation}\notag
\xymatrix{ Rf_*\Lam_{C_v} \otimes^{\L} \T_r(1)_{\scX} \ar[rr]^{\id \otimes^{\L} \, \res^f \quad}
 \ar[dd]_{\gys_f \otimes^{\L} \, \id } && Rf_*\Lam_{C_v} \otimes^{\L} \, Rf_*\T_r(1)_{C_v}
 \ar[d]^{{\text{product}}} \\ &&  Rf_* \T_r(1)_{C_v} \ar[d]^{\gys_f} \\
\T_r(N)_{\scX}\,[2N] \otimes^{\L} \T_r(1)_{\scX} \ar[rr]^-{\text{product}} && \T_r(d)_{\scX}\,[2N].  }
\end{equation}
This completes the proof of Lemma \ref{lem:CD} and Theorem \ref{thm:surj}.
\qed


\medskip
\begin{thebibliography}{CTSia}
\bibitem[Be]{Be}
       Beilinson, A. A.:
       Height pairings between algebraic cycles.
       In: Manin, Yu.\ I. (ed.) {\it $K$-theory, Arithmetic and Geometry},
       (Lecture Notes in Math.\ 1289), pp.\ 1--27,
       Berlin, Springer, 1987
\bibitem[Bl]{Bl}
       Bloch, S.:
       On the Chow groups of certain rational surfaces.
       Ann.\ Sci.\ Ec.\ Norm.\ Sup.\ (4) {\bf 14}, 41--59 (1981) 
\bibitem[BK]{bk}
       Bloch, S., Kato, K.:
       $p$-adic \'etale cohomology.
       Inst.\ Hautes \'Etudes Sci.\ Publ.\ Math.\ {\bf 63}, 107--152 (1986)
\bibitem[CT1]{CT1}
       Colliot-Th\'el\`ene, J.-L.:
       Hilbert's Theorem 90 for $K_2$, with application to the Chow groups of rational surfaces.
       Invent.\ Math.\ {\bf 71}, 1--20 (1983).
\bibitem[CT2]{CT0}
       Colliot-Th\'{e}l\`{e}ne, J.-L.:
       Cycles alg\'{e}briques de torsion et
       $K$-th\'{e}orie alg\'{e}brique.
       In: Ballico, E.\ (ed.) {\it Arithmetic Algebraic Geometry, Trento 1991}
       (Lecture notes in Math.\ {\bf 1553}), pp.\ 1--49, Berlin, 1993
\bibitem[CT3]{CT2}
       Colliot-Th\'el\`ene, J.-L.:
       L'arithm\'etique du groupe de Chow des z\'ero-cycles.
       J. Th\'eor.\ Nombres Bordeaux {\bf 7}, 51--73 (1995)
\bibitem[CT4]{CTb}
      Colliot-Th\'el\`ene, J.-L.:
      Groupe de Chow des z\'ero-cycles sur les vari\'et\'es $p$-adiques (d'apr\`es S. Saito, K. Sato et al.).
      S\'eminaire N. Bourbaki 2009--2010, n$^o$ 1012
\bibitem[CT5]{ct:letter}
      Colliot-Th\'el\`ene, J.-L.:
      Letter to the first author.
\bibitem[CTOP]{CTOP}
      Colliot-Th\'el\`ene, J.-L., Ojanguren, M., Parimala, R.:
      Quadratic forms over fraction fields of two-dimensional henselian rings and Brauer groups of related schemes.
      In: Parimala, R. (ed.) {\it Algebra, Arithmetic and Geometry, Mumbai 2000},
      (Tata Inst.\ Fund.\ Res.\ Stud.\ Math.\ 16, Vol.\ I), pp.\ 185--217,
      Tata Inst.\ Fund.\ Res., Bombay, 2002.
\bibitem[CTS]{cts}
       Colliot-Th\'{e}l\`{e}ne, J.-L., Saito, S.:
       Z\'ero-cycles sur les vari\'et\'es $p$-adiques et groupe de Brauer.
       Internat.\ Math.\ Res.\ Notices 1996, 151--160
\bibitem[CTSn]{CTSn}
       Colliot-Th\'{e}l\`{e}ne, J.-L., Sansuc, J.-J.:
       On the Chow groups of certain rational surfaces: a sequel to a paper of S. Bloch.
       Duke Math.\ J. {\bf 48}, 421--447 (1981)
\bibitem[CTSS]{CTSS}
       Colliot-Th\'{e}l\`{e}ne, J.-L., Sansuc, J.-J., Soul\'{e}, C.:
       Torsion dans le groupe de Chow de codimension deux.
       Duke Math.\ J. {\bf 50}, 763--801 (1983)
\bibitem[Da]{Da}
       Dalawat, C. S.:
       Le groupe de Chow d'une surface rationnelle sur un corps local.
       Compositio Math.\ {\bf 141}, 344--358 (2005)
\bibitem[dJ]{dJ}
       de Jong, A. J.:
       Smoothness, semi-stability and alterations.
       Inst.\ Hautes \'Etudes Sci.\ Publ.\ Math.\ {\bf 83}, 51--93 (1996)
\bibitem[F]{Fu}
       Fujiwara, K.:
       A proof of the absolute purity conjecture (after Gabber).
       In: Usui, S., Green, M., Illusie, L., Kato, K.,
       Looijenga, E., Mukai, S., Saito, S. (eds.)
       {\it Algebraic Geometry, Azumino, 2001},
       (Adv.\ Stud.\ in Pure Math.\ 36), pp.\ 153--184,
       Tokyo, Math.\ Soc.\ Japan, 2002
\bibitem[Ga]{gabber}
       Gabber, O.:
       Some theorems on Azumaya algebras.
       In: Kervaire, M., Ojanguren, M. (eds.)
       {\it Groupe de Brauer S\'eminaire, Les Plans-sur-Bex, 1980},
       (Lecture Notes in Math.\ 844), pp.\ 129--209,
       Berlin, Springer, 1981
\bibitem[Gr]{Gr1}
       Gros, M.:
       Classes de Chern et classes des cycles en cohomologie logarithmique.
       Bull.\ Soc.\ Math.\ France M\'{e}moire N${}^o$ 21, 1985
\bibitem[G]{gr}
       Grothendieck, A.:
       Le groupe de Brauer.
       In: {\it Dix Expos\'es sur la Cohomologie des Sch\'emas},
       pp.\ 46--188,
       Amsterdam, North-Holland, 1968
\bibitem[H]{Hy}
       Hyodo, O.:
       A note on $p$-adic etale cohomology in the semi-stable reduction case.
       Invent.\ Math.\ {\bf 91} 543--557 (1988)
\bibitem[J]{Ja1}
       Jannsen, U.:
       Continuous \'etale cohomology.
       Math.\ Ann.\ {\bf 280}, 207--245 (1987)
\bibitem[K]{kk:cft}
        Kato, K.:
        A generalization of local class field theory by using $K$-groups II.
        J. Fac.\ Sci.\ Univ.\ of Tokyo Sec.\ IA {\bf 27}, 603--683 (1980)
\bibitem[KS]{KS}
        Kato, K., Saito, S.:
        Global class field theory of arithmetic schemes.
        In: {\it Applications of algebraic K-theory to algebraic geometry and number theory, Boulder, 1983}
        (Contemporary Math.\ 55 Part I), pp.\ 255--331, Providence, Amer.\ Math.\ Soc., 1986
\bibitem[Kh]{kahn}
        Kahn, B.:
        Descente galoisienne et $K_2$ des corps de nombres. $K$-Theory {\bf 7}, 55--100 (1993)
\bibitem[KL]{KL}
        Katz, N., Lang, S.:
        Finiteness theorem for higher geometric class field theory.
        Enseign.\ Math.\ {\bf 27}, 285--319 (1981)
\bibitem[Ke]{Ke}
        Kerz, M.: Id\`eles in higher dimension. preprint, {\tt http://arxiv.org/abs/0907.5337}
\bibitem[L1]{Li2}
        Lichtenbaum, S.:
        Duality theorems for curves over $p$-adic fields.
        Invent.\ Math.\ {\bf 7}, 120--136 (1969)
\bibitem[L2]{Li}
       Lichtenbaum, S.:
       Values of zeta functions at non-negative integers.
       In: Jager, H. (ed.) {\it Number Theory, Noordwijkerhout, 1983},
       (Lecture Notes in Math.\ 1068), pp.\ 127-138,
       Berlin, Springer, 1984
\bibitem[M1]{M}
       Manin, Yu. I.:
       Le groupe de Brauer-Grothendieck en g\'eom\'etrie diophantinne.
       In: {\it Actes du Congr\'es International des Math\'ematicients,
       Nice 1970}, Tome 1, pp.\ 401--411, Paris, Gauthier-Villars, 1971
\bibitem[M2]{M2}
       Manin, Yu. I.:
       {\it Cubic forms. Algebra, geometry, arithmetic}. Translated by Hazewinkel, M., 2nd ed.
       (North-Holland Math.\ Library 4), North-Holland, Amsterdam, 1986  
\bibitem[Ma]{mat}
       Mattuck, A.:
       Abelian varieties over $\p$-adic ground fields.
       Ann.\ of Math.\ {\bf 62}, 92--119 (1955)
\bibitem[MS]{MS}
       Merkur'ev, A. S., Suslin, A. A.:
       $K$-cohomology of Severi-Brauer varieties
       and the norm residue homomorphism.
       Math.\ USSR Izv.\ {\bf 21}, 307--341 (1983)
\bibitem[Mi]{Mi2}
       Milne, J. S.:
       Values of zeta functions of varieties over finite fields.
       Amer.\ J. Math.\ {\bf 108}, 297--360 (1986)
\bibitem[PS]{PS}
       Parimala R., Suresh, V.:
       Zero-cycles on quadratic fibrations:
       Finiteness theorems and the cycle map.
       Invent.\ Math.\ {\bf 122}, 83--117 (1995)
\bibitem[RZ]{RZ}
    Rapoport, M., Zink, T.:
    \"{U}ber die lokale Zetafunktion von Shimuravariet\"{a}ten.
  Monodromiefiltration und verschwindende Zyklen in ungleicher Charakteristik.
    Invent.\ Math.\ {\bf 68}, 21--101 (1982)
\bibitem[S1]{ss:localring}
    Saito, S.:
    Arithmetic on two-dimensional local rings. Invent.\ Math.\ {\bf 85}, 379--414 (1986)
\bibitem[S2]{ss:cycle}
    Saito, S.:
    On the cycle map for torsion algebraic cycles of codimension two. Invent.\ Math.\ {\bf 106}, 443--460 (1991)
\bibitem[SS1]{SS}
    Saito, S., Sato, K.:
    A finiteness theorem for zero-cycles over $p$-adic fields. Ann.\ of Math.\ {\bf 172}, 1593--1639 (2010)
\bibitem[SS2]{SaSa}
    Saito, S., Sato, K.:
    A $p$-adic regulator map and finiteness results for arithmetic schemes.
    Documenta Math.\ Extra Volume: Andrei A. Suslin's Sixtieth Birthday, 525--594 (2010)
\bibitem[Sa1]{sato1}
    Sato, K.:
    Logarithmic Hodge-Witt sheaves on normal crossing varieties.
    Math.\ Z.\ {\bf 257}, 707--743 (2007)
\bibitem[Sa2]{Sat2}
    Sato, K.:
    $p$-adic \'etale Tate twists and arithmetic duality.
    (with an appendix by Hagihara, K.),
    Ann.\ Sci.\ \'Ecole Norm.\ Sup.\ (4) {\bf 40}, 519--588 (2007)
\bibitem[Sch]{Sch}
    Schneider, P.:
    $p$-adic point of motives.
    In: Jannsen, U. (ed.) {\it Motives}, (Proc.\ Symp.\ Pure Math.\ 55-II), pp.\ 225--249,
    Providence, Amer.\ Math.\ Soc., 1994
\bibitem[Se1]{Se}
    Serre, J.-P.:
    {\it Cohomologie Galoisienne.} 5${}^{e}$ \'ed., (Lecture Notes in Math.\ 5),
    Berlin, Springer, 1992
\bibitem[Se2]{Se2}
    Serre, J.-P.:
    Zeta and $L$-functions.
    In: Schilling, O. F. G. (ed.) {\it Arithmetical Algebraic Geometry}, pp.\ 82--92,
    New York, Harper and Row, 1965
\bibitem[Sh]{Sh}
    Shiho, A.:
    On logarithmic Hodge-Witt cohomology of regular schemes.
    J. Math.\ Sci.\ Univ.\ Tokyo {\bf 14}, 567--635 (2007)
\bibitem[Sp]{Spa}
    Spaltenstein, N.:
    Resolutions of unbounded complexes.
    Compositio Math.\ {\bf 65}, 121--154 (1988)
\bibitem[Sr]{srinivas}
    Srinivas, V.:
    {\it Algebraic $K$-Theory. 2nd ed.},
    (Progr.\ Math.\ 90),
    Boston, Birkh\"{a}user, 1996
\bibitem[Th]{Th}
    Thomason, R. W.:
    Absolute cohomological purity.
    Bull.\ Soc.\ Math.\ France {\bf 112}, 397--406 (1984)
\bibitem[Y1]{Y1}
    Yamazaki, T.:
    Formal Chow groups, $p$-divisible groups and syntomic cohomology. Duke Math.\ J. {\bf 102}, 359--390 (2000)
\bibitem[Y2]{Y2}
    Yamazaki, T.:
    On Chow and Brauer groups of a product of Mumford curves. Math.\ Ann.\ {\bf 333}, 549--567 (2005)
\bibitem[SGA2]{SGA2}
   Grothendieck, A.:
   {\it Cohomologie locale des faisceaux coh\'erents et th\'eor\`emes de Lefschetz locaux et globaux.}
   Documents Math\'ematiques (Paris)  4, Paris, Soc.\ Math.\ France, 2005
\bibitem[SGA4]{sga4}
    Grothendieck, A., Artin, M., Verdier, J.-L.,
    with Deligne, P., Saint-Donat, B.:
    {\it Th\'{e}orie des Topos et Cohomologie \'{E}tale des Sch\'{e}mas, Tome 3}.
    (Lecture Notes in Math.\ 305),
    Berlin, Springer, 1973
\end{thebibliography}
\end{document}